\documentclass{article}

\usepackage{amsmath,amssymb,bm}
\usepackage{cite,url,color}

\usepackage[top=2.5cm,bottom=2.5cm,left=2cm,right=2cm]{geometry}

\usepackage{longtable}

\newtheorem{definition}{Definition}
\newtheorem{notation}{Notation}
\newtheorem{theorem}{Theorem}
\newtheorem{lemma}{Lemma}
\newtheorem{corollary}{Corollary}

\newtheorem{example}{Example}
\newtheorem{remark}{Remark}
\newtheorem{problem}{Open Problem \#}
\newenvironment{proof}[1][Proof]{\par\noindent\textit{#1}: }{\hfill$\blacksquare$\vskip 0.5\baselineskip}
\newenvironment{solution}[1][Solution]{\par\noindent\textit{#1}: }{\hfill$\square$\vskip 0.5\baselineskip}

\newcommand\fproblem[1]{%
\begin{center}
\fbox{\parbox{0.95\textwidth}{%
\begin{problem}
#1
\end{problem}}}
\end{center}}

\usepackage[normalem]{ulem}

\usepackage{arydshln}

\begin{document}

\title{Permutation Polynomials modulo $m$}
\author{Shujun Li\\\url{http://www.hooklee.com}}
\date{\today}

\maketitle

\begin{abstract}
This paper mainly studies problems about so called ``permutation
polynomials modulo $m$", polynomials with integer
coefficients\footnote{In this paper, we always call them ``integer
polynomial" in short. Note that sometime another name ``integral
polynomial" is used \cite[Sec. 7.2]{HW:NumberTheory1979}. However,
we prefer to ``integer polynomial" to avoid confusion with the word
``integral" as an adjective (see \cite{IntegerPolynomial}).} that
can induce bijections over $\mathbb{Z}_m=\{0,\cdots,m-1\}$. The
necessary and sufficient conditions of permutation polynomials are
given, and the number of all permutation polynomials of given degree
and the number induced bijections are estimated. A method is
proposed to determine all equivalent polynomials from the induced
polynomial function, which can be used to determine all equivalent
polynomials that induce a given bijection. A few problems have not
been solved yet in this paper and left for open study.

\textit{Note: After finishing the first draft, we noticed that some
results obtained in this paper can be proved in other ways (see
Remark \ref{remark:RepeatingWork}). In this case, this work gives
different and independent proofs of related results.}
\end{abstract}

\tableofcontents

\setlength{\arrayrulewidth}{0.6pt} \setlength{\doublerulesep}{0pt}
\renewcommand\arraystretch{1.25}

\section{Introduction}

Integer polynomials that can induce bijections over finite fields,
namely permutation polynomials, are firstly studied in algebra
community \cite{LauschNobauer:AlgebraPoly1973,
LidlMullen:PermPolyI:AMM1988, LidlMullen:PermPolyII:AMM1990,
RMT:DicksonPoly1993, Mullen:PermPolySurvey:FFaTA1995, SunWan:FF1987,
Lidl:FF1997}. Permutation polynomials have been used in cryptography
and coding \cite{PPinRSA-Crypt83, PP-Ciphers-EuroCrypt84,
MI-PP-PKC-IEEEISIT83, Cryptanalysis-PP-PKC-EuroCrypt84,
Cade-PP-PKC-SIAM-ALA85, BreakingCade-Crypto86,
ModifyingCade-Crypto86, SunWan:FF1987, Mullen:PP-NFSR:IEEETIT1989,
PP2Coding2005}. This paper studies permutation polynomials modulo an
integer, i.e., permutation polynomials over integer rings
\cite{Mullen:PolyFun-mod:AMH1984, RMT:DicksonPoly1993,
Rivest:PPmod2w:FFTA2001, Sun:PPmodCoding:IEEETIT2005}. Such
permutation polynomials have also been used in cryptography and
coding recently, such as in the RC6 block cipher \cite{RC6} a simple
permutation polynomial $f(x)=x(2x+1)$ modulo $2^d$ is used.

Assume $f(x)=a_nx^n+\cdots+a_1x+a_0$ is a polynomial with integer
coefficients of degree $n\geq 1$ modulo $m$, where $a_n\not\equiv
0\pmod m$. It is possible that $f(x)$ forms a bijection over
$\mathbb{Z}_m=\{0,\cdots,m-1\}$, i.e., $\forall
x_1,x_2\in\mathbb{Z}$ and $x_1\not\equiv x_2\pmod m$,
$f(x_1)\not\equiv f(x_2)\pmod m$. In other words, it is true that
$f(\mathbb{Z})=\mathbb{Z}_m$, or a complete system of residues
modulo $m$ is permuted by the polynomial $f$. The most common
permutation polynomial modulo $m$ is $f(x)=x$. In addition, as a
special case, Fermat's little theorem also gives one of the simplest
permutation polynomials modulo a prime $p$: $f(x)=x^p$, which
satisfies $\forall x\in\mathbb{Z}$, $f(x)\equiv x\pmod p$.

It is well-known that many problems on permutation polynomials over
finite fields are still open \cite{LidlMullen:PermPolyI:AMM1988,
LidlMullen:PermPolyII:AMM1990, Mullen:PermPolySurvey:FFaTA1995}.
Similarly, there are a few work on permutation polynomials modulo
integers, both in number theory \cite{HW:NumberTheory1979,
Eynden:NumberTheory1987, NZM:NumberTheory1991,
Rosen:NumberTheory1993, Rose:NumberTheory1994,
Pan:ConciseNumberTheory1998} and algebra communities. In this paper,
we try to find answers to the following questions on permutation
polynomials modulo $m$, where $m$ may be a prime, a prime power, or
a general composite.
\begin{itemize}
\item
What are necessary and sufficient conditions of permutation
polynomials modulo $m$?

\item
What is the number of distinct permutation polynomials of degree
$\leq n$ modulo $m$ and what is the number of distinct bijections
induced from these polynomials?

\item
Is there a practical way to enumerate all permutation polynomials of
degree $n$ modulo $m$ given a bijection (or a partial bijection)
over $\{0,\cdots,m-1\}$?
\end{itemize}

At present, the first two questions have been almost solved, but the
last one has not been solved at all when the degree of the
polynomials $n\geq p$ modulo $p^d$.

This paper is organized as follows. In next section, we first give
some preliminary definitions and lemma. Sec. \ref{section:composite}
discusses the case of composite moduli, and concludes that
permutation polynomials modulo a composite can be studied via
permutation polynomials modulo each prime power. In Sec.
\ref{section:prime}, we discuss the case of prime moduli and give
some limited results. Then, in Sec. \ref{section:PrimePowers}, the
general cases modulo $p^d$ ($d\geq 1$) are studied and some useful
results are obtained. Two open problems are raised in Sec.
\ref{section:PrimePowers} for future study.

\section{Preliminaries}
\label{section:preliminaries}

This section lists a number of definitions and notations used
throughout in this paper. Some preliminary lemmas are also given to
simplify the discussions in this paper. I try to keep the
definitions, notations and lemmas as simple as possible. Please feel
free to contact me if you have some idea of making them even
simpler, more elegant, more beautiful, and/or more rigorous in
mathematics.

\subsection{Some Simple Lemmas on Congruences}

The following lemmas will be extensively cited in this paper without
explicit citations.

\begin{lemma}\label{lemma:dvi-mod}
If $a\mid b$ and $a\mid m$, then $a\mid (b\bmod m)$.
\end{lemma}
\begin{proof}
From $a|b$, $\exists k_1\in\mathbb{Z}$, $b=ak_1$. From $a|m$,
$\exists k_2\in\mathbb{Z}$, $m=ak_2$. Assume $x=(b\bmod m)$, then
$\exists k_3\in\mathbb{Z}$, $x=mk_3+b=ak_2k_3+ak_1=a(k_1+k_2k_3)$.
So, $a\mid x$, which proves this lemma.
\end{proof}

\begin{lemma}\label{lemma:mod}
If $m\mid m'$, $(a\bmod m')\equiv a \pmod m$.
\end{lemma}
\begin{proof}
From $m|m'$, $\exists k_1\in\mathbb{Z}$, $m'=mk_1$. Assume $x=
a\bmod m'$, then $\exists k_2\in\mathbb{Z}$, $x=m'k_2+a=mk_1k_2+a$.
Then, $(a\bmod m')\equiv mk_1k_2+a\equiv a\pmod m$. This lemma is
proved.
\end{proof}

\begin{lemma}
If $a\equiv 0\pmod{m_1}$ and $b\equiv 0\pmod{m_2}$, then $ab\equiv
0\pmod{m_1m_2}$.
\end{lemma}

\begin{lemma}[Theorem 2.2 in \cite{NZM:NumberTheory1991}]\label{lemma:polynomial-mod}
Assume $f(x)=a_nx^n+\cdots+a_1x+a_0$ is an integer polynomial. If
$x_1\equiv x_2\pmod m$, then $f(x_1)\equiv f(x_2)\pmod m$.
\end{lemma}

\begin{lemma}
Assume $\bm{A}$ is an $n\times n$ matrix, $\bm{X}$ is a vector of
$n$ unknown integers, and $\bm{B}$ is a vector of $n$ integers. If
$|\bm{A}|$ is relatively prime to $m$, i.e., $\gcd(|\bm{A}|,m)=1$,
then $\bm{A}\bm{X}\equiv\bm{B}\pmod m$ has a unique set of
incongruent solutions
$\bm{X}\equiv\overline{\Delta}(\mathrm{adj}(\bm{A}))\bm{B}\pmod m$,
where $\overline{\Delta}$ is an inverse of $\Delta=|\bm{A}|$ modulo
$m$ and $\mathrm{adj}(\bm{A})$ is the adjoint of $\bm{A}$.
\end{lemma}
\begin{proof}
This lemma is a direct result of Theorem 3.18 in
\cite{Rosen:NumberTheory1993} (see pages 151 and 152).
\end{proof}

\subsection[Polynomial Congruences Modulo $m$]%
{Polynomial Congruences Modulo $\bm{m}$}

The following definition is from Chap. VII of
\cite{HW:NumberTheory1979} and related concepts are slightly
extended.

\begin{definition}
Given two integer polynomials of degree $n$:
$f(x)=a_nx^n+\cdots+a_1x+a_0$ and $g(x)=b_nx^n+\cdots+b_1x+b_0$, if
$\forall i=0\sim n$, $a_i\equiv b_i\pmod m$, we say \uline{$f(x)$ is
congruent to $g(x)$ modulo $m$}, or $f(x)$ and $g(x)$ are
\uline{congruent (polynomials) modulo $m$}, which is denoted by
$f(x)\equiv g(x)\pmod m$. On the other hand, if $\exists
i\in\{1,\cdots,n\}$, such that $a_i\not\equiv b_i\pmod m$, we say
$f(x)$ and $g(x)$ are \uline{incongruent (polynomials) modulo $m$},
denoted by $f(x)\not\equiv g(x)\pmod m$.
\end{definition}

\begin{definition}
A \uline{polynomial congruence (residue) class modulo $m$} is a set
of all polynomials congruent to each other modulo $m$.
\end{definition}

\begin{definition}
A set of polynomials of degree $n$ modulo $m$ is a \uline{complete
system of polynomial residues of degree $n$ modulo $m$}, if for
every polynomial of degree $n$ modulo $m$ there is one and only one
congruent polynomial in this set.
\end{definition}

\begin{lemma}
The following set of polynomials is a complete system of polynomial
residues of degree $n$ modulo $m$:
\[
\mathbb{F}[x]=\left\{f(x)=a_nx^n+\cdots+a_1x+a_0\left|a_n\in\{1,\cdots,m-1\},a_{n-1},\cdots,a_0\in\{0,\cdots,m-1\}\right.\right\}.
\]
\end{lemma}
\begin{proof}
Assume $f(x)=a_nx^n+\cdots+a_1x+a_0$ is a polynomial of degree $n$
modulo $m$. Choose $a_i^*=(a_i\bmod m)\in\{0,\cdots,m-1\}$ ($i=0\sim
n$), then $f^*(x)=a_n^*x^n+\cdots+a_1^*x+a_0^*\in\mathbb{F}$ is
congruent to $f(x)$. Assume that another polynomial
$g(x)=b_nx^n+\cdots+b_1x+b_0\in\mathbb{F}$ is also congruent to
$f(x)$. Then, $\forall i=0\sim n$, $b_i\equiv a_i^*\pmod m$. Since
$\{0,\cdots,m-1\}$ is a complete set of residues modulo $m$,
$b_i=a_i^*$. This means that $g(x)=f^*(x)$. This completes the proof
of this lemma.
\end{proof}

\begin{definition}
A set of polynomials of degree $\leq n$ modulo $m$ is a
\uline{complete system of polynomial residues of degree $\leq n$
modulo $m$}, if for every polynomial of degree $\leq n$ modulo $m$
there is one and only one congruence polynomial.
\end{definition}

\begin{lemma}
The following set of polynomials is a complete system of polynomial
residues of degree $n$ modulo $m$:
\[
\mathbb{F}[x]=\left\{f(x)=a_nx^n+\cdots+a_1x+a_0\left|a_n,a_{n-1},\cdots,a_0\in\{0,\cdots,m-1\}\right.\right\}.
\]
\end{lemma}
\begin{proof}
The proof is similar to the above lemma.
\end{proof}

\subsection[Polynomial Functions Modulo $m$]{Polynomial Functions
Modulo $\bm{m}$}

\begin{definition}
If a function over $\{0,\cdots,m-1\}$ can be represented by a
polynomial modulo $m$, we say this function is \uline{polynomial
modulo $m$}.
\end{definition}

\begin{lemma}\label{lemma:any-fun-poly-p}
Assume $p$ is a prime. Then, any function over $\{0,\cdots,p-1\}$ is
polynomial modulo $p$.
\end{lemma}
\begin{proof}
Assume $f(x)=a_nx^n+\cdots+a_1x+a_0$ is a polynomial of degree
$n\geq p-1$ modulo $p$. Given a function
$F:\{0,\cdots,p-1\}\to\{0,\cdots,p-1\}$, one has the following
system of congruences:
\[
\left[\begin{matrix}%
1 & 0 & 0^2 & \cdots & 0^{p-1}\\
1 & 1 & 1^2 & \cdots & 1^{p-1}\\
1 & 2 & 2^2 & \cdots & 2^{p-1}\\
\vdots & \vdots & \vdots & \ddots & \vdots\\
1 & p-1 & (p-1)^2 & \cdots & (p-1)^{p-1}
\end{matrix}\right]
\left[\begin{matrix}%
a_0\\
a_1\\
a_2\\
\vdots\\
a_{p-1}
\end{matrix}\right]\equiv
\left[\begin{matrix}%
F(0)\\
F(1)-\sum_{i=p}^na_i\\
F(2)-\sum_{i=p}^n2^ia_i\\
\vdots\\
F(p-1)-\sum_{i=p}^n(p-1)^ia_i
\end{matrix}\right]\pmod{p^d}.
\]
Since the matrix at the left side is a Vondermonde matrix, one can
see its determinant is relatively prime to $p$. So, for each
combination of $a_p,\cdots,a_n$, there is a unique set of
incongruent solutions of $a_0,\cdots,a_{p-1}$. Thus this lemma is
proved.
\end{proof}

\subsection[Equivalent Polynomials Modulo $m$]{Equivalent Polynomials Modulo $\bm{m}$}

The concept of equivalent polynomial modulo $m$ is used to describe
incongruent but equivalent (for any integer) polynomials modulo $m$.
Note that some researchers call them ``residually congruent
polynomials modulo $m$" \cite{Kempner:PolyResidue:TAMS1921a,
Kempner:PolyResidue:TAMS1921b}.

\begin{definition}
Two integer polynomials $f(x)$ and $g(x)$ are \uline{equivalent
(polynomials) modulo $m$} if $\forall x\in\mathbb{Z}$, $f(x)\equiv
g(x)\pmod m$.
\end{definition}
Note that two equivalent polynomials modulo $m$ may not be congruent
modulo $p$, and may have distinct degrees. As a typical example,
when $p$ is a prime, $f(x)=x^p$ and $g(x)=x$ are equivalent
polynomials modulo $p$.

\begin{lemma}\label{lemma:equivalent-poly-degree1}
Two polynomials of degree 1 modulo $m$, $f(x)=a_1x+a_0$ and
$g(x)=b_1x+b_0$, are equivalent polynomials modulo $m$ if and only
if $f(x)\equiv g(x)\pmod m$, i.e., $a_1\equiv b_1\pmod m$ and
$a_0\equiv b_0\pmod m$.
\end{lemma}
\begin{proof}
The ``if" part is obvious from the definition of equivalent
polynomials modulo $m$, so we focus on the ``only if" part. Since
$f(x)$ and $g(x)$ are equivalent polynomials modulo $m$, then
$\forall x\in\{0,\cdots,m-1\}$,
$f(x)-g(x)=(a_1-b_1)x+(a_0-b_0)\equiv 0\pmod m$. Choosing $x\equiv
0\pmod m$, one has $a_0\equiv b_0\pmod m$. Then, choosing $x\equiv
1\pmod m$, one has $a_1\equiv b_1\pmod m$. Thus this lemma is
proved.
\end{proof}

\begin{lemma}\label{lemma:equivalent-poly-a0}
Two polynomials, $f(x)=a_{n_1}x^{n_1}+\cdots+a_0$ and
$g(x)=b_{n_2}x^{n_2}+\cdots+b_0$, are equivalent polynomials modulo
$m$, then $a_0\equiv b_0\pmod m$.
\end{lemma}
\begin{proof}
Choosing $x=0$, one has $f(x)-g(x)=a_0-b_0\equiv 0\pmod m$. This
lemma is proved.
\end{proof}
\begin{corollary}
Two polynomials, $f(x)=a_nx^n+\cdots+a_2x^2+a_0$ and
$g(x)=a_nx^n+\cdots+a_2x^2+b_0$, are equivalent polynomials modulo
$m$ if and only if $a_0\equiv b_0\pmod m$.
\end{corollary}

\begin{lemma}\label{lemma:equivalent-poly-pd}
Assume $p$ is a prime and $d\geq 1$. Two polynomials,
$f(x)=a_{p-1}x^{p-1}+\cdots+a_0$ and
$g(x)=b_{p-1}x^{p-1}+\cdots+b_0$, are equivalent polynomials modulo
$p^d$ if and only if $f(x)\equiv g(x)\pmod{p^d}$.
\end{lemma}
\begin{proof}
The ``if" part is obvious true, from the definition of equivalent
polynomials modulo $p^d$. So, we focus on the ``only if" part only.
From $f(x)-g(x)\equiv 0\pmod{p^d}$, choosing $x=0\sim p-1$, one can
get the following system of congruences in the matrix form
$\bm{A}\bm{X}_{a-b}\equiv \bm{B}\pmod{p^d}$:
\begin{equation}
\left[\begin{matrix}%
1 & 0 & 0^2 & \cdots & 0^{p-1}\\
1 & 1 & 1^2 & \cdots & 1^{p-1}\\
1 & 2 & 2^2 & \cdots & 2^{p-1}\\
\vdots & \vdots & \vdots & \ddots & \vdots\\
1 & p-1 & (p-1)^2 & \cdots & (p-1)^{p-1}
\end{matrix}\right]
\left[\begin{matrix}%
a_0-b_0\\
a_1-b_1\\
a_2-b_2\\
\vdots\\
a_{p-1}-b_{p-1}
\end{matrix}\right]\equiv
\left[\begin{matrix}%
f(0)-g(0)\\
f(1)-g(1)\\
f(2)-g(2)\\
\vdots\\
f(p-1)-g(p-1)
\end{matrix}\right]\equiv
\left[\begin{matrix}%
0\\
0\\
0\\
\vdots\\
0
\end{matrix}\right]\pmod{p^d}.\label{equation:equivalent-poly-pd}
\end{equation}
Since $\bm{A}$ is a Vandermonde sub-matrix, one can get
$|\bm{A}|=\prod_{0\leq i<j\leq p-1}(j-i)$
\cite[\S4.4]{Zhang:MatrixTheory1999}. From $p$ is a prime and
$1\leq(j-i)\leq p-1$, one has $\gcd(|\bm{A}|,p^d)=1$. Thus, the
above system of congruences has a unique set of incongruent
solutions. So, $\forall i=0\sim p-1$, one has $a_i\equiv
b_i\pmod{p^d}$. This completes the proof of this lemma.
\end{proof}
Note that in the above lemma $f(x)$ and $g(x)$ may be polynomials of
degree less than $p-1$ modulo $p^d$. In this case, the matrix at the
left side of the system of congruences may have a smaller size, but
its determinant is still relatively prime to $p^d$.

\begin{corollary}
Assume $p$ is a prime. Two polynomials, $f(x)=a_nx^n+\cdots+a_0$ and
$g(x)=b_nx^n+\cdots+b_0$, are equivalent polynomials modulo $p$ if
and only if $(f(x)\bmod(x^p-x))\equiv(g(x)\bmod(x^p-x))\pmod p$.
\end{corollary}
\begin{proof}
This corollary is a direct result of the above lemma and Fermat's
Little Theorem.
\end{proof}

\subsection[Permutation Polynomials Modulo $m$]
{Permutation Polynomials Modulo $\bm{m}$}

\begin{definition}
Assume $f(x)=a_nx^n+\cdots+a_1x+a_0$ is a polynomial of degree
$n\geq 1$ modulo $m$, where $a_n\not\equiv 0\pmod m$. If
$f(x)=((a_nx^n+\cdots+a_1x+a_0)\bmod m)$ forms a bijection
$F:\{0,\cdots,m-1\}\to\{0,\cdots,m-1\}$, we say that $f(x)$ is a
\uline{permutation polynomial modulo $m$}, or $f(x)$ is
\uline{permutation modulo $m$}. The bijection $F$ is called the
\uline{induced bijection of the polynomial $f(x)$ modulo $m$}.
\end{definition}

\begin{definition}
If two permutation polynomials are equivalent modulo $m$, we say
they are \uline{equivalent permutation polynomials modulo $m$}. It
is obvious that equivalent permutation polynomials modulo $m$ induce
the same bijection over $\{0,\cdots,m-1\}$.
\end{definition}

\begin{lemma}\label{lemma:bijective}
A polynomial $f(x)$ is a permutation polynomial modulo $m$ if and
only if $g(x)=af(x)+b$ is a permutation polynomial modulo $m$, where
$\gcd(a,m)=1$ and $b\in\mathbb{Z}$.
\end{lemma}
\begin{proof}
This lemma is a direst result of Theorem 3.6 in
\cite{Rosen:NumberTheory1993} on a complete system of residues
modulo $m$.
\end{proof}
\begin{lemma}\label{lemma:bijective-poly-a0}
Two polynomials, $f(x)=a_{n_1}x^{n_1}+\cdots+a_1x+a_0$ and
$f(x)=b_{n_2}x^{n_2}+\cdots+b_1x+b_0$, is equivalent permutation
polynomials if and only if $f^*(x)=a_{n_1}x^{n_1}+\cdots+a_1x$ and
$g^*(x)=b_{n_2}x^{n_2}+\cdots+b_1x$ are equivalent permutation
polynomials modulo $m$ and $a_0\equiv b_0\pmod m$.
\end{lemma}
\begin{proof}
This lemma is a direct result of Lemmas
\ref{lemma:equivalent-poly-a0} and \ref{lemma:bijective}.
\end{proof}
From the above two lemmas, we can only study permutation polynomials
in the form $f(x)=a_nx^n+\cdots+a_1x$.

\begin{lemma}
If $f(x)\equiv g(x)\pmod m$ and $f(x)$ is a permutation polynomial
modulo $m$, then $g(x)$ is an equivalent permutation polynomial of
$f(x)$ modulo $m$.
\end{lemma}
\begin{proof}
Since $f(x)\equiv g(x)\pmod m$, $\forall a\in\mathbb{Z}$,
$f(a)\equiv g(a)\pmod m\Rightarrow f(a)\bmod m=g(a)\bmod m$, i.e.,
$F(a)=G(a)$. So $g(x)$ generates the same bijection as $f(x)$. This
completes the proof.
\end{proof}

\begin{theorem}\label{theorem:bijection-degree1}
The polynomial $f(x)=a_1x+a_0$ is a permutation polynomial modulo
$m$ if and only if $\gcd(a_1,m)=1$.
\end{theorem}
\begin{proof}
Assume $G=\langle 1\rangle$ is a cyclic group of order $m$. From
Theorem 3.24 in \cite{Gilbert:Algebra2005} (or Theorem 2 in \S 2.3
of \cite{Hu:Algebra1999}), $a_1=1^{a_1}$ is a generator of $G$ if
and only if $\gcd(a_1,m)=1$. Note that in group $G$ the binary
operator is defined as addition modulo $m$. It is obvious that $a_1$
is a generator of $G$ if and only if $f^*(x)=a_1x$ is a permutation
polynomial modulo $m$. Then from Lemma \ref{lemma:bijective}, this
theorem is proved.
\end{proof}

\begin{corollary}\label{corollary:bijective-poly-degree1}
The number of congruence classes of permutation polynomials
$f(x)=a_1x+a_0$ of degree 1 modulo $m$ is $\phi(m)m$. The number of
bijections induced from these permutation polynomials is also
$\phi(m)m$.
\end{corollary}
\begin{proof}
From Theorem \ref{theorem:bijection-degree1}, $a_1$ should satisfy
$\gcd(a_1,m)=1$, but $a_0$ can be any integer, so The number of
congruence classes of permutation polynomials of degree 1 modulo $m$
is $\phi(m)m$. From Lemma \ref{lemma:equivalent-poly-degree1}, the
$\phi(m)m$ permutation polynomials are not equivalent to each other,
so they induce $\phi(m)m$ distinct bijections.
\end{proof}


\begin{definition}
Given a bijection $F:\mathbb{A}\to\mathbb{A}$. If for a set
$\mathbb{B}\subseteq\mathbb{A}$, $F(\mathbb{B})=\mathbb{B}$, then
$F_\mathbb{B}=\{(a,b)|a,b\in\mathbb{B}\}\subseteq F$ is a bijection
over $\mathbb{B}$, and we say the bijection
$F_\mathbb{B}:\mathbb{B}\to\mathbb{B}$ is a \uline{sub-bijection} of
$F$, and $F$ is a \uline{super-bijection} of $F_\mathbb{B}$.
\end{definition}

\begin{lemma}\label{lemma:sub-bijection}
Given a bijection $F:\mathbb{A}\to\mathbb{A}$. If
$F_\mathbb{B}:\mathbb{B}\to\mathbb{B}$ is a sub-bijection of $F$,
then
$F_{\mathbb{A}\backslash\mathbb{B}}:\mathbb{A}\backslash\mathbb{B}\to\mathbb{A}\backslash\mathbb{B}$
is also a sub-bijection of $F$.
\end{lemma}
\begin{proof}
Assume that $\exists a\in\mathbb{A}\backslash\mathbb{B}$,
$F(a)\in\mathbb{B}$. Since $F_\mathbb{B}$ is a bijection over
$\mathbb{B}$, then $F(a)$ has one and only one preimage in
$\mathbb{B}$. However, it is obvious that $a\not\in\mathbb{B}$ is
also the preimage of $F(a)$. We get a contradiction. So, $\forall
x\in\mathbb{A}\backslash\mathbb{B}$,
$F(x)\in\mathbb{A}\backslash\mathbb{B}$. This means that
$F_{\mathbb{A}\backslash\mathbb{B}}$ is a sub-bijection of $F$ over
$\mathbb{A}\backslash\mathbb{B}$.
\end{proof}

\subsection[Null Polynomials modulo $m$]{Null Polynomials modulo $\bm{m}$}

This concept was introduced in \cite{Li:NullPoly2005}, and also
studied by others without a special name
\cite{Kempner:PolyResidue:TAMS1921a, Kempner:PolyResidue:TAMS1921b}.
Here, we just give the definition and some simple lemmas on null
polynomials modulo $m$. For more advanced results, see
\cite{Li:NullPoly2005}.

\begin{definition}
A polynomial $f(x)$ of degree $n\geq 0$ modulo $m$ is a \uline{null
polynomial of degree $n$ modulo $m$}, if $\forall x\in\mathbb{Z}$,
$f(x)\equiv 0\pmod m$. Specially, $f(x)=0$ is a trivial null
polynomial of degree 0 modulo $m$.
\end{definition}

\begin{lemma}
If $f(x)=a_nx^n+\cdots+a_1x+a_0$ is a null polynomial modulo $m$,
then $a_0\equiv 0\pmod m$.
\end{lemma}

\begin{lemma}
Given any null polynomial $f(x)$ modulo $m$, $af(x)$ will still be a
null polynomial modulo $m$, where $a$ is an arbitrary integer.
\end{lemma}
\begin{lemma}
A polynomial $f(x)$ is a null polynomial modulo $m$, if and only
$af(x)$ is a null polynomial modulo $m$, where $\gcd(a,m)=1$.
\end{lemma}
\begin{lemma}\label{lemma:null-poly-transitivity}
If $f(x)$ is a null polynomial modulo $m$ and $a\mid m$, then $f(x)$
is still a null polynomial modulo $a$.
\end{lemma}
The most frequently used form of the above lemma is as follows: if
$f(x)$ is a null polynomial modulo $p^d$, then $f(x)$ is still a
null polynomial modulo $p^i$ for any integer $i\leq d$.

\begin{lemma}\label{lemma:Equ-Poly-Null-Poly}
Two polynomials, $f_1(x)$ and $f_2(x)$, are equivalent polynomials
modulo $m$ if and only if $f_1(x)-f_2(x)$ is a null polynomial
modulo $m$.
\end{lemma}

\begin{definition}
Denote the least integer $n\geq 1$ such that there exists a null
polynomial of degree $n$ modulo $m$ by $\omega_0(m)$ and call it
\uline{the least null-polynomial degree modulo $m$}. Denote the
least integer $n\geq 1$ such that there exists a \textbf{monic} null
polynomial of degree $n$ modulo $m$ by $\omega_1(m)$ and call it
\uline{the least monic null-polynomial degree modulo $m$}. A (monic)
null polynomial of degree $\omega_0(m)$ or $\omega_1(m)$ is called
\uline{a least-degree (monic) null polynomial modulo $m$}. 
\end{definition}

\begin{lemma}
Every polynomial of degree $\geq\omega_1(m)$ modulo $m$ has one
equivalent polynomial of degree $\leq\omega_1(m)-1$ modulo $m$.
\end{lemma}

\subsection[Circular Shift of an Integer Set Modulo $m$ (New)]%
{Circular Shift of an Integer Set Modulo $\bm{m}$ (New)}

\begin{definition}
The \uline{$k$-th circular shift} of an integer set $\mathbb{A}$
\uline{modulo $m$} is defined by $\{x|x=(y+k)\bmod
m,y\in\mathbb{A}\}$ and denoted by $[(\mathbb{A}+k)\bmod m]$ in this
paper.
\end{definition}

\begin{lemma}\label{lemma:circular-shift}
If $\mathbb{A}_1,\cdots,\mathbb{A}_k$ is a partition of an integer
set $\mathbb{A}=\{0,\cdots,m-1\}$, then $\forall a\in\mathbb{Z}$,
$[(\mathbb{A}_1+a)\bmod m],\cdots,[(\mathbb{A}_k+a)\bmod m]$ is
still a partition of $\mathbb{A}$.
\end{lemma}
\begin{proof}
From $\mathbb{A}_1,\cdots,\mathbb{A}_k$ is a partition of
$\mathbb{A}$, one has
$\bigcup_{i=1}^k\mathbb{A}_i=\{0,\cdots,m-1\}$. Then,
$\bigcup_{i=1}^k[(\mathbb{A}_i+a)\bmod
m]=\bigcup_{i=1}^k\{x|x=(y+a)\bmod
m,y\in\mathbb{A}_i\}=\left\{x|x=(y+a)\bmod
m,y\in\bigcup_{i=1}^k\mathbb{A}_i\right\}=\{x|x=(y+a)\bmod
m,y\in\mathbb{A}\}$. From $\mathbb{A}$ is a complete system of
residues modulo $m$, $\forall x_0\in\mathbb{A}$ and $\forall
a\in\mathbb{Z}$, $\exists y_0\in\mathbb{A}$, $x_0-a\equiv y_0\pmod
m\Leftrightarrow x_0\equiv y_0+a\pmod m$. Since $x_0\in\mathbb{A}$,
i.e., $0\leq x_0\leq m-1$, one has $x_0=(y_0+a)\bmod m$. This means
$\forall x_0\in\bigcup_{i=1}^k[(\mathbb{A}_i+a)\bmod
m]=\{x|x=(y+a)\bmod m,y\in\mathbb{A}\}$. So
$\bigcup_{i=1}^k[(\mathbb{A}_i+a)\bmod m]=\mathbb{A}$.

On the other hand, since $\mathbb{A}_1,\cdots,\mathbb{A}_k$ is a
partition of $\mathbb{A}$, $\forall y_1\in\mathbb{A}_i$ and $\forall
y_2\in\mathbb{A}_j$ ($i\neq j$), one has $y_1\neq y_2$. Considering
$\mathbb{A}$ is a complete system of residues modulo $m$, one
immediately gets $y_1\not\equiv y_2\pmod m
\Leftrightarrow(y_1+a)\not\equiv (y_2+a)\pmod m \Leftrightarrow
(y_1+a)\bmod m\neq (y_2+a)\bmod m$. This means
$\mathbb{A}_i\cap\mathbb{A}_j=\varnothing$.

The above two results proves this lemma.
\end{proof}

\subsection[Base-$p$ Resolution (New)]{Base-$\bm{p}$ Resolution (New)}

\begin{definition}
The \uline{base-$p$ resolution} of an integer $a$ is an integer
$i\geq 0$ such that $p^i\parallel a$, i.e., $p^i\mid a$ but
$p^{i+1}\nmid a$. Specially, define the base-$p$ resolution of 0 as
$+\infty$. When $p=2$, the base-$p$ resolution is also called the
\uline{binary resolution}.
\end{definition}

\begin{notation}
The set of all integers of base-$p$ resolution $i$ is denoted by
$\mathbb{Z}(i|_p)$. The set of all elements in $\mathbb{A}$ of
base-$p$ resolution $i$ is denoted by $\mathbb{A}(i|_p)$. The set of
all elements in $\mathbb{A}$ of base-$p$ resolution $i\geq a$ is
denoted by $\mathbb{A}(\geq a|_p)$; similarly, we can define
$\mathbb{A}(\leq a|_p)$ and $\mathbb{A}(\neq a|_p)$. The set of all
elements in $\mathbb{A}$ of base-$p$ resolution $a\leq i\leq b$ is
denoted by $\mathbb{A}(a\mapsto b|_p)$. The set of all elements in
$\mathbb{A}$ of base-$p$ resolution $i\in\{i_1,\cdots,i_k\}$ is
denoted by $\mathbb{A}(i_1,\cdots,i_k|_p)$. In the above notations,
the subscription ``$p$" denotes the base (radix) of the resolution.
\end{notation}
Apparently, $a\in\mathbb{Z}(\geq i|_p)\Leftrightarrow
p^i|a\Leftrightarrow a\equiv 0\pmod{p^i}$.

\begin{definition}
The \uline{base-$p$ multi-resolution partition} of the integer set
$\mathbb{Z}$ is a collection of the following sets:
$\{\mathbb{Z}(i|_p)\}_{i=0}^{+\infty}$. The \uline{base-$p$
multi-resolution partition} of an integer set $\mathbb{A}$ is a
collection of the following $(k_{\max}-k_{\min}+2)$ sets:
$\mathbb{A}(k_{\min}|_p),\cdots,\mathbb{A}(k_{\max}|_p)$ and
$\mathbb{A}(+\infty|_p)=\{0\}$, where $k_{\min}$ and $k_{\max}$
denote the minimal and the maximal base-$p$ resolution of all
non-zero integers in $\mathbb{A}$.
\end{definition}
Specially, the base-$p$ multi-resolution partition of
$\mathbb{A}=\{0,\cdots,p^d-1\}$ is a collection of the following
$d+1$ sets: $\mathbb{A}(0|_p),\cdots,\mathbb{A}(d-1|_p)$ and
$\mathbb{A}(+\infty|_p)=\{0\}$.

\begin{definition}
Assume the base-$p$ resolution of an integer $a\neq 0$ is $i\geq 0$,
then the \uline{base-$p$ representation of $a$} is a sequence of $i$
integers $a_0,\cdots,a_i$, such that $a=\sum_{j=0}^ia_jp^j$ and
$a_j\in\{0,\cdots,p-1\}$. It is denoted by $a=(a_i\cdots a_0)_p$.
Specially, the base-$p$ representation of 0 is $(0)_p$. The $j$-th
integer in the base-$p$ representation of $a$ is called \uline{the
$j$-th base-$p$ digit} or \uline{the $j$-th digit of base $p$} or
\uline{the $j$-th digit} in short if the base is well defined in the
context.
\end{definition}
It is obvious that the base-$p$ resolution of an integer is unique
and $a_j=\lfloor a/p^j\rfloor\bmod p$.

\subsection{Determinants of Some Special Matrices}

\begin{lemma}\label{lemma:det-powers}
Assume $m\geq 1$. Given a $2m\times 2m$ matrix
$\bm{A}=\left[\begin{matrix}\bm{A}_1\\\bm{A}_2\end{matrix}\right]$,
where $\bm{A}_1=[X_j^{i-1}]_{1\leq j\leq m \atop 1\leq i\leq 2m}$
and $\bm{A}_2=[iX_j^{i-1}]_{1\leq j\leq m \atop 1\leq i\leq 2m}$,
i.e.,
\[
\bm{A}=\left[\begin{array}{cccc:ccc}%
1 & X_1 & \cdots & X_1^{m-1} & X_1^m & \cdots & X_1^{2m-1}\\
1 & X_2 & \cdots & X_2^{m-1} & X_2^m & \cdots & X_2^{2m-1}\\
\vdots & \vdots & \ddots & \vdots & \vdots & \ddots & \vdots\\
1 & X_m & \cdots & X_m^{m-1} & X_m^m & \cdots & X_m^{2m-1}\\
\hdashline
1 & 2X_1 & \cdots & mX_1^{m-1} & (m+1)X_1^m & \cdots & 2mX_1^{2m-1}\\
1 & 2X_2 & \cdots & mX_2^{m-1} & (m+1)X_2^m &
\cdots & 2mX_2^{2m-1}\\
\vdots & \vdots & \ddots & \vdots & \vdots & \ddots & \vdots\\
1 & 2X_m & \cdots & mX_m^{m-1} & (m+1)X_m^m & \cdots & 2mX_m^{2m-1}
\end{array}\right].
\]
Then, $|\bm{A}|=(-1)^{\frac{m(m-1)}{2}}\prod_{j=1}^mX_j\prod_{1\leq
i<j\leq m}(X_j-X_i)^4$.
\end{lemma}
\begin{proof}
A proof can be found in \cite{Li:Determinants2005}, or in
\cite{ADC1999} (as a special case of Theorem 20).
\end{proof}
\begin{corollary}\label{corollary:det-powers}
Assume $m\geq 1$. Given a $2m\times 2m$ matrix
$\bm{A}=\left[\begin{matrix}\bm{A}_1\\\bm{A}_2\end{matrix}\right]$,
where $\bm{A}_1=[X_i^{j+1}]_{1\leq i\leq m \atop 1\leq j\leq 2m}$
and $\bm{A}_2=[(j+1)X_i^j]_{1\leq i\leq m \atop 1\leq j\leq 2m}$.
Then,
$|\bm{A}|=(-1)^{\frac{m(m-1)}{2}}\prod_{i=1}^mX_i^4\prod_{1\leq
i<j\leq m}(X_j-X_i)^4$.
\end{corollary}

\begin{lemma}\label{lemma:det-binom-powers}
Assume $m\geq 1,n\geq l\geq 1$ and $\bm{A}$ is a block-wise
$ml\times ml$ matrix as follows:
\[
\bm{A}=\left[\begin{matrix}%
\bm{A}_1\\
\bm{A}_2\\
\vdots\\
\bm{A}_m
\end{matrix}\right],
\]
where for $i=1\sim m$,
\[
\bm{A}_i=\left[\binom{n+j-1}{k-1}X_i^{j-1}\right]_{1\leq
j\leq ml \atop 1\leq k\leq l}=\left[\begin{matrix}%
\binom{n}{0} & \binom{n+1}{0}X_i & \cdots &
\binom{n+(ml-1)}{0}X_i^{ml-1}\\
\binom{n}{1} & \binom{n+1}{1}X_i & \cdots &
\binom{n+(ml-1)}{1}X_i^{ml-1}\\
\vdots & \vdots & \ddots & \vdots\\
\binom{n}{l-1} & \binom{n+1}{l-1}X_i & \cdots &
\binom{n+(ml-1)}{l-1}X_i^{ml-1}
\end{matrix}\right]_{ml\times l}.
\]
Then, $|\bm{A}|=\prod_{i=1}^mX_i^{\frac{l(l-1)}{2}}\prod_{1\leq
i<j\leq m}(X_j-X_i)^{l^2}$.
\end{lemma}

\section{Permutation Polynomials modulo $\bm{m=p_1^{d_1}\cdots p_r^{d_r}}$}
\label{section:composite}

The theorems given in this section says that we can focus our study
on permutation polynomials modulo prime and prime powers.

\begin{theorem}\label{theorem:bijective-poly-composite}
Assume $p_1$, $\cdots$, $p_r$ are $r$ distinct prime numbers and
$d_1$, $\cdots$, $d_r\geq 1$. A polynomial $f(x)$ is a permutation
polynomial modulo $m=\prod_{i=1}^rp_i^{d_i}$, if and only if
$\forall i=1\sim r$, $f(x)$ is a permutation polynomial modulo
$p_i^{d_i}$.
\end{theorem}
\begin{proof}
To simplify the following proof, $\forall i=1\sim r$, assume
$P_i=p_i^{d_i}$ and $\overline{P_i}=m/P_i$. Since $p_1,\cdots,p_r$
are all primes, it is obvious that $\gcd(P_i,\overline{P_i})=1$. In
addition, assume $\mathbb{M}=\{0,\cdots,m-1\}$, and $\forall i=1\sim
r$, $\mathbb{M}_i=\{0,\cdots,P_i-1\}$ and
$\mathbb{M}_i^*=\{0,\cdots,\overline{P_i}-1\}$.

First, $\forall i=1\sim r$, let us prove the ``only if" part.
$\forall a\in\mathbb{M}_i^*$, assume
$\mathbb{A}=\{x|x\in\mathbb{M},x\equiv
a\pmod{\overline{P_i}}\}=\{\overline{P_i}y+a|y\in\mathbb{M}_i\}$.
Since $f(x)$ is a permutation polynomial modulo
$m=P_i\overline{P_i}$, $\forall y_1,y_2\in\mathbb{A}$ and $y_1\neq
y_2$, one has $f(y_1)\not\equiv f(y_2)\pmod{P_i}$ or
$f(y_1)\not\equiv f(y_2)\pmod{\overline{P_i}}$, otherwise it
conflicts with part 3) of Theorem 2.3 in \cite{NZM:NumberTheory1991}
(i.e., Property IX in \S15 of \cite{Pan:ConciseNumberTheory1998}).
Since $f(y_1)\equiv f(y_2)\equiv a\pmod{\overline{P_i}}$, one
immediately knows $f(y_1)\not\equiv f(y_2)\pmod{P_i}$. From Lemma
\ref{theorem:bijection-degree1}, since $\gcd(\overline{P_i},P_i)=1$,
$g(y)=\overline{P_i}y+a$ is a permutation polynomial modulo $P_i$.
This means that $\mathbb{A}$ is a complete system of residues modulo
$P_i$. This leads to the result that $f(x)$ is a permutation
polynomial modulo $P_i$.

Next, we prove the ``if" part. Given $r$ integers as follows:
$a_1\in\mathbb{M}_1$, $\cdots$, $a_r\in\mathbb{M}_r$, construct the
system of $r$ simultaneous congruences, $i=1\sim r:f(x)\equiv
a_i\pmod{P_i}$. From the Chinese Remainder Theorem, there is exactly
one solution of $f(x)$ in each complete system of residues modulo
$m$. Since $f(x)$ is a permutation polynomial modulo each $P_i$, we
can construct $m=\prod_{i=1}^rp_i^{d_i}$ systems of $r$ simultaneous
congruences, and get $m$ distinct solutions of $f(x)$ in each
complete system of residues modulo $m$. Considering there are only
$m$ elements in each complete system of residues modulo $m$, one can
immediately deduce that $f(x)$ is also a permutation polynomial
modulo $m$.
\end{proof}

\begin{theorem}
Assume $p_1$, $\cdots$, $p_r$ are $r$ distinct prime numbers, $d_1$,
$\cdots$, $d_r\geq 1$ and $m=\prod_{i=1}^rp_i^{d_i}$. If $f_1(x)$,
$\cdots$, $f_r(x)$ are permutation polynomials modulo $p_1^{d_1}$,
$\cdots$, $p_r^{d_r}$, respectively, then there exists one and only
one permutation polynomial $f(x)$ modulo $m$ in each complete system
of polynomial residues modulo $m$, such that $f(x)\equiv
f_i(x)\pmod{p_i^{d_i}}$ holds for $i\in\{1,\cdots,r\}$.
\end{theorem}
\begin{proof}
Applying the Chinese remainder theorem on each coefficients of the
polynomials, one can immediately prove this theorem.
\end{proof}

\section{Permutation Polynomials modulo $\bm{p}$}
\label{section:prime}

It is natural to connect Fermat's Little Theorem with permutation
polynomials, since this theorem actually says that there always
exists a permutation polynomial $f(x)=x^p$ of degree $p$ modulo a
prime $p$ such that $\forall x\in\mathbb{Z}$, $f(x)\equiv x\pmod p$.
However, the original Fermat's Little Theorem say nothing about how
many permutation polynomials there are and how to calculate other
permutation polynomials (if any). We have an enhanced version to
answer this question.

\begin{theorem}\label{theorem:Fermat-enhanced}
Assume $p$ is a prime. There exist $(p-1)p!$ congruence classes of
permutation polynomials of degree $p$ modulo $p$. For each given
bijection over $\{0,\cdots,p-1\}$, there exist $p-1$ congruence
classes of permutation polynomial of degree $p$ modulo $p$.
Specially, there exists a permutation polynomial $f(x)=x^p$, such
that $\forall x\in\mathbb{Z}$, $f(x)\equiv x\pmod p$.
\end{theorem}
\begin{proof}
Assume $f(x)=a_px_p+a_{p-1}x^{p-1}+\cdots+a_1x+a_0$, where
$a_p\not\equiv 0\pmod p$. Choosing $x=0,\cdots,p-1$, respectively,
one can get the following $p$ congruences modulo $p$.
\begin{eqnarray*}
a_p\cdot 0^p+a_{p-1}\cdot 0^{p-1}+\cdots a_1\cdot 0+a_0 & \equiv &
f(0)\pmod p\\
a_p\cdot 1^p+a_{p-1}\cdot 1^{p-1}+\cdots a_1\cdot 1+a_0 & \equiv &
f(1)\pmod p\\
& \vdots\\
a_p\cdot (p-1)^p+a_{p-1}\cdot (p-1)^{p-1}+\cdots a_1\cdot(p-1)+a_0 &
\equiv & f(p-1)\pmod p
\end{eqnarray*}
Fixing $a_p$, rewrite the above system of congruences as the
following matrix form $\bm{A}\bm{X}_a\equiv \bm{B}\pmod p$.
\begin{equation}
\left[\begin{matrix}%
1 & 0 & 0 & \cdots & 0\\
1 & 1 & 1 & \cdots & 1\\
1 & 2 & 2^2 & \cdots & 2^p\\
\vdots & \vdots & \vdots & \ddots & \vdots\\
1 & p-1 & (p-1)^2 & \cdots & (p-1)^p
\end{matrix}\right]
\left[\begin{matrix}%
a_0\\
a_1\\
a_2\\
\vdots\\
a_{p-1}
\end{matrix}\right]\equiv
\left[\begin{matrix}%
f(0)\\
f(1)-a_p\\
f(2)-2^pa_p\\
\vdots\\
f(p-1)-(p-1)^pa^p
\end{matrix}\right]\pmod p\label{equation:Fermat-enhanced}
\end{equation}
Apparently, $\bm{A}$ is a Vandermonde matrix, so its determinant can
be calculated as $|\bm{A}|=\prod_{0\leq i<j\leq p-1}(j-i)$
\cite[\S4.4]{Zhang:MatrixTheory1999}. Since $p$ is a prime and
$0\leq(j-i)\leq p-1$, one has $\gcd(|\bm{A}|,p)=1$. Thus, the above
system of congruence has a unique (i.e., one and only one) solution
modulo $p$, for each combination of the values of
$f(0),\cdots,f(p-1)$ and $a_p$. For each possible value of $a_p$,
the number of all possible combinations of the values of
$f(0),\cdots,f(p-1)$ is $p!$. Since $a_p$ has $p-1$ congruence
classes modulo $p$, one immediately deduces that there exists
$(p-1)p!$ congruence classes of permutation polynomials of degree
$p$ modulo $p$. For each bijection over $\{0,\cdots,p-1\}$, i.e.,
for each combination of the values of $f(0),\cdots,f(p-1)$, there
are $p-1$ distinct congruence classes of permutation polynomials of
degree $p$ modulo $p$, each of which corresponds to one possible
value of $a_p$ modulo $p$.

When $a_p=1$, choosing $f(i)=i$ for $i=0\sim p-1$, one can get a
special solution: $a_0\equiv a_1\equiv \cdots\equiv a_{p-1}\equiv
0\pmod p$. This leads to $f(x)=x^p\equiv x\pmod p$, which is the
permutation polynomial of degree $p$ modulo $p$ as mentioned in the
Fermat's little theorem.
\end{proof}

From the above theorem, one can get some more results on the number
of permutation polynomials and induced bijections modulo $p$.

\begin{notation}
Assume $p$ is a prime. Denote the number of distinct permutation
polynomials and the number of all distinct polynomials in a complete
system of polynomial resides of degree $\leq n$ modulo $p$ by
$N_{pp}(\leq n,p)$ and $N_p(\leq n,p)$ respectively. Here, the
subscript ``pp" means ``permutation polynomial" and ``p" denotes
``polynomial". Similar functions will be defined later.
\end{notation}

\begin{corollary}\label{corollary:count-bijective-poly-p}
Assume $p$ is a prime. The following is true: when $n\geq p-1$,
$\dfrac{N_{pp}(\leq n,p)}{N_p(\leq n,p)}=\dfrac{(p-1)!}{p^{p-1}}$.
\end{corollary}
\begin{proof}
Recall the proof of Theorem \ref{theorem:Fermat-enhanced}, when
$n\geq p-1$, changing the degree of the polynomial from $p$ to $n$
and moving $a_p,\cdots,a_n$ to the right side, Eq.
(\ref{equation:Fermat-enhanced}) has a unique set of incongruent
solutions to the values of $f(0),\cdots,f(p-1),a_p,\cdots,a_n$.
Since $f(0),\cdots,f(p-1)$ forms a complete permutation modulo $p$,
one immediately has $N_{pp}(\leq n)/N_p(\leq
n)=p!/p^p=(p-1)!/p^{p-1}$.
\end{proof}

\begin{corollary}\label{corollary:count-bijection-p}
Assume $p$ is a prime and $n\geq p-1$. The number of bijections
induced from permutation polynomials of degree $\leq n$ modulo $p$
is $p!$.
\end{corollary}
\begin{proof}
This corollary can be proved in the same way as the above corollary,
due to the fact that each permutation of
$f(0),\cdots,f(p-1),a_p,\cdots,a_n$ corresponds to a unique set of
incongruent solutions to $a_0,\cdots,a_{p-1}$.
\end{proof}

\begin{corollary}[A special case of Lemma \ref{lemma:equivalent-poly-pd}]
Assume $p$ is a prime. Two permutation polynomials of degree $\leq
p-1$ modulo $p$, $f_1(x)$ and $f_2(x)$, are equivalent if and only
if they are congruence polynomials modulo $p$, i.e., $f_1(x)\equiv
f_2(x)\pmod p$.
\end{corollary}
\begin{proof}
The ``if" part is obvious. Let us see the ``only if" part. From the
above two corollaries, the number of permutation polynomials of
degree $\leq p-1$ modulo $p$ and the number of bijections induced
from these polynomials are both $p!$. This immediately leads to the
fact that any two equivalent permutation polynomials are congruent
polynomials, otherwise the number of bijections will be less than
$p!$. Thus, this corollary is true.
\end{proof}

\begin{corollary}
Assume $p$ is a prime and $f(x)=a_nx^n+\cdots+a_1x+a_0$ is a
permutation polynomial of degree $n\geq p$ modulo $p$. Then $f(x)$
has exactly $p^{n-p}$ equivalent polynomials of degree $\leq n$
modulo $p$ (including itself).
\end{corollary}
\begin{proof}
This corollary can be proved in a similar way to the above
corollaries.
\end{proof}

\section{Permutation Polynomials modulo $\bm{p^d}$ ($\bm{d\geq 1}$)}
\label{section:PrimePowers}

\subsection{Hierarchy Theorem}

This theorem shows the hierarchical structure of the bijection
induced from a permutation polynomial modulo $p^d$.

\begin{theorem}[Hierarchy Theorem]\label{theorem:hierarchy-pd}
Assume $p$ is a prime and $f(x)=a_nx^n+\cdots+a_1x$ is a permutation
polynomial of degree $n$ modulo $m=p^d$ and
$\mathbb{A}=\{0,\cdots,p^d-1\}$. The following results are true.
\begin{enumerate}
\item
The induced bijection $F:\mathbb{A}\to\mathbb{A}$ is composed of two
sub-bijections, $F_1:\mathbb{A}(\geq 1|_p)\to\mathbb{A}(\geq 1|_p)$
and $F_0:\mathbb{A}\backslash\mathbb{A}(\geq
1|_p)\to\mathbb{A}\backslash\mathbb{A}(\geq 1|_p)$.

\item
The sub-bijection $F_0$ is composed of $p-1$ sub-bijections,
$\forall i=1\sim p-1$, $F_{0,i}:[(\mathbb{A}(\geq 1|_p)+i)\bmod
m]\to\left[\left(\mathbb{A}(\geq
1|_p)+\sum_{k=1}^ni^ka_k\right)\bmod m\right]$.

\item
When $d\geq 2$, each of the $p$ sub-bijections, $F_1$ and
$F_{0,1},\cdots,F_{0,p-1}$, corresponds to a permutation polynomial
of degree $\leq d-1$ modulo $p^{d-1}$ in the following form:
$f^*(z)=\sum_{i=d-1}^1b_ip^{i-1}z^i=b_{d-1}p^{d-2}z^{d-1}+\cdots+b_2pz+b_1z$.

\item
When $d\geq 2$, $F_1$ has is composed of $d$ sub-bijections:
$i=1\sim d-1$, $F_{1,i}:\mathbb{A}(i|_p)\to\mathbb{A}(i|_p)$, and
$F_{1,d}:\mathbb{A}(+\infty|_p)\to\mathbb{A}(+\infty|_p)$. Each of
$F_{0,1},\cdots,F_{0,p-1}$ is also composed of $d$ sub-bijections of
this kind.

\item
When $d\geq 1$, $\forall i,j\in\{0,\cdots,p-1\}$ and $i\neq j$,
$\sum_{k=1}^na_k(j^k-i^k)=a_1(j-i)+a_2(j^2-i^2)+\cdots+a_n(j^n-i^n)\not\equiv
0\pmod p$.

\item
When $d\geq 2$, $\forall i\in\{0,\cdots,p-1\}$,
$\sum_{k=1}^nki^{k-1}a_k=a_1+2\cdot i^1\cdot a_2+\cdots+n\cdot
i^{n-1}\cdot a_n\not\equiv 0\pmod p$.
\end{enumerate}
\end{theorem}
\begin{proof}
We prove all the results one by one. Note that $\mathbb{A}(\geq
1|_p)$, $[(\mathbb{A}(\geq 1|_p)+1)\bmod m],\cdots,[(\mathbb{A}(\geq
1|_p)+p-1)\bmod m]$ forms a partition of $\mathbb{A}$.

1. $\forall x\in\mathbb{A}(\geq 1|_p)$, then $p\mid x$. Since $x\mid
f(x)$, so $p\mid f(x)$. This means that $f(x)$ forms a sub-bijection
$F_1:\mathbb{A}(\geq 1|_p)\to\mathbb{A}(\geq 1|_p)$. From Lemma
\ref{lemma:sub-bijection}, there exists another sub-bijection
$F_0:\mathbb{A}\backslash\mathbb{A}(\geq
1|_p)\to\mathbb{A}\backslash\mathbb{A}(\geq 1|_p)$.

2 \& 5. $\forall i=1\sim p-1$ and $\forall x\in[(\mathbb{A}(\geq
1|_p)+i)\bmod p]$, one has $x\equiv i\mod p$. Assume $x=y+i$, where
$y\in\mathbb{A}(\geq 1|_p)$, and one can get a new polynomial as
follows:
\[
f_{0,i}^*(y)=f(y+i)=a_n(y+i)^n+\cdots+a_1(y+i)=f_{0,i}^{**}(y)+\sum_{k=1}^ni^ka_k,
\]
where
$f_{0,i}^{**}(y)=\sum_{l=n}^1\left(\sum_{k=l}^n\binom{k}{l}i^{k-l}a_ky^l\right)$.
Applying the first result on $f_{0,i}^{**}(y)$, one knows it forms a
sub-bijection over $\mathbb{A}(\geq 1|_p)$. So, $f(x)$ forms a
sub-bijection $F_{0,i}:[(\mathbb{A}(\geq 1|_p)+i)\bmod
p]\to\left[\left(\mathbb{A}(\geq
1|_p)+\sum_{k=1}^ni^ka_k\right)\bmod m\right]$. This proves the 2nd
result of this theorem. Since $f(x)$ is a permutation polynomial
modulo $p^d$, $\forall i\neq j$, the ranges of
$F_{0,1},\cdots,F_{0,p-1}$ should form a partition of
$\mathbb{A}\backslash\mathbb{A}(\geq 1|_p)$. This means that
$\forall i,j\in\{0,\cdots,p-1\}$ and $i\neq j$,
$\sum_{k=1}^ni^ka_k\not\equiv \sum_{k=1}^nj^ka_k\pmod p$. This leads
to the 5th result of this theorem.

3. For $F_1$, since $p\mid x$, let us assume $x=pz$, where
$z\in\{0,\cdots,p^{d-1}-1\}$. Substitute $x=pz$ into $f(x)$, we have
another polynomial
$f_1(z)=f(pz)=a_n(pz)^n+\cdots+a_1(pz)=pf_1^*(z)$, where
$f_1^*(z)=a_np^{n-1}z^n+\cdots+a_1z$. Apparently, over
$\mathbb{A}(\geq |_p)$, $f(x)$ is uniquely determined by the
polynomial $f_1^*(z)$ modulo $p^{d-1}$. This means that $f_1^*(z)$
is a permutation polynomial modulo $p^{d-1}$. When $n\geq d$,
$p^{d-1}\mid a_np^{n-1}z^n+\cdots+a_dp^{d-1}z^d$, so $f_1^*(z)\equiv
a_{d-1}p^{d-2}z^{d-1}+\cdots+a_1z\pmod{p^{d-1}}$. As a result, the
degree of $f_1^*(z)$ modulo $p^{d-1}$ is always not greater than
$d-1$. For $F_{0,1},\cdots,F_{0,p-1}$, applying the same analysis on
$f_{0,1}^{**}(y),\cdots,f_{0,p-1}^{**}(y)$, one can get a similar
result.

4 \& 6. When $x\in\mathbb{A}(+\infty|_p)=\{0\}$, $f(x)=f(0)=0$, so
there exists a sub-bijection
$F_{1,d}:\mathbb{A}(+\infty|_p)\to\mathbb{A}(+\infty|_p)$. $\forall
i=1\sim d-1$ and $\forall x\in\mathbb{A}(i|_p)$, then $\exists
k_1,k_2\in\mathbb{Z}$ and $k_2\not\equiv 0\pmod p$, such that
$x=p^i(k_1p+k_2)$. Then,
$x^2=p^{2i}(k_1p+k_2)^2=p^{i+1}p^{i-1}(k_1p+k_2)^2$, so $p^{i+1}\mid
x^2$. As a result, $f(x)\equiv a_1x=a_1p^i(k_1p+k_2)\pmod{p^{i+1}}$.
Assume $a_1\equiv 0\pmod p$, one has $f(x)\equiv 0\pmod{p^{i+1}}$.
This means that $f(x)\in\mathbb{A}(\geq i+1|_p)\backslash\{0\}$.
However, since $p\geq 2$, the cardinality of $\mathbb{A}(\geq
i+1|_p)\backslash\{0\}$ is always smaller than the cardinality of
$\mathbb{A}(i|_p)$, which conflicts with the fact that $f(x)$ is a
permutation polynomial modulo $p^d$. So one immediately has
$a_1\not\equiv 0\pmod p$ and $f(x)\in\mathbb{A}(i|_p)$, i.e., $f(x)$
forms a sub-bijection over $\mathbb{A}(i|_p)$. For $i=1\sim p-1$,
applying the same analysis for $f_{0,i}^{**}(y)$, we can get similar
results: $\sum_{k=1}^nki^{k-1}a_k\not\equiv 0\pmod p$ and
$f(x)_{0,i}^{**}(y)\in\mathbb{A}(i|_p)$. Thus the 4th and the 6th
results have been proved.
\end{proof}

\subsection{Necessary and Sufficient Conditions for Permutation Polynomials}

\begin{theorem}\label{theorem:bijective-poly-NSC-pd}
Assume $p$ is a prime and $d\geq 2$. The polynomial
$f(x)=a_nx^n+\cdots+a_1x$ is a permutation polynomial modulo $p^d$
if and only if the following two conditions are true simultaneously:
\begin{enumerate}
\item
$f(x)$ is a permutation polynomial modulo $p$, i.e., $\forall
i,j\in\{0,\cdots,p-1\}$ and $i\neq j$,
$f(j)-f(i)=\sum_{k=1}^na_k(j^k-i^k)=a_1(j-i)+a_2(j^2-i^2)+\cdots+a_n(j^n-i^n)\not\equiv
0\pmod p$.

\item
$\forall i\in\{0,\cdots,p-1\}$,
$\sum_{k=1}^nki^{k-1}a_k=a_1+2ia_2+\cdots+ni^{n-1}a_n\not\equiv
0\pmod p$.
\end{enumerate}
\end{theorem}
\begin{proof}
The ``only if" part of this theorem has been proved in Theorem
\ref{theorem:hierarchy-pd}, so we only focus on the ``if" part. Let
us use mathematical induction on $d$ to prove this part.

1) When $d=2$, consider the $p$ sub-bijections, $F_1$,
$F_{0,1},\cdots,F_{0,p-1}$, separately.

When $x\in\mathbb{A}(\geq 1|_p)$, assume $x=pz$, where
$z\in\{0,\cdots,p-1\}$, so $f(x)=a_nx^n+\cdots+a_2x^2+a_1x\equiv
p\cdot a_1z\pmod{p^2}$. Apparently, $f(x)$ is uniquely determined by
the polynomial $f_1^*(z)=a_1z$ modulo $p$. Choosing $i=0$, the
second necessary and sufficient condition becomes $a_1\not\equiv
0\pmod p$, which means $\gcd(a_1,p)=1$. Then, from Lemma
\ref{theorem:bijection-degree1}, $f_1^*(z)=a_1z$ forms a bijection
over $\{0,\cdots,p-1\}$. This means $f(x)$ forms a bijection $F_1$
over $\mathbb{A}(\geq 1|_p)$.

$\forall i=1\sim p-1$ and $\forall x\in[(\mathbb{A}(\geq
1|_p)+i)\bmod p]$, assume $x=y+i$, where $y\in\mathbb{A}(\geq
1|_p)$. Substitute $x=y+i$ into $f(x)$, one has
$f_{0,i}^*(y)=f_{0,i}^{**}(y)+\sum_{k=1}^ni^ka_k$, where
$f_{0,i}^{**}(y)=\sum_{l=n}^1\left(\sum_{k=l}^n\binom{k}{l}i^{k-l}a_ky^l\right)$.
Assume $y=pz$, where $z\in\{0,\cdots,p-1\}$, due to the same reason
in the case of $x\in\mathbb{A}(\geq 1|_p)$, $f_{0,i}^{**}(y)\equiv
\sum_{k=1}^nki^{k-1}a_ky\equiv p\sum_{k=1}^nki^{k-1}a_kz\pmod{p^2}$.
The second necessary and sufficient condition ensure that
$\gcd\left(\sum_{k=1}^nki^{k-1}a_k,p\right)=1$, so
$\sum_{k=1}^nki^{k-1}a_kz$ forms a bijection over $\{0,\cdots,p-1\}$
and thus $f_{0,i}^{**}(y)$ forms a bijection over $\mathbb{A}(\geq
1|_p)$. This further leads to the fact that $f(x)$ forms a bijection
$F_{0,i}:[(\mathbb{A}(\geq 1|_p)+i)\bmod
p]\to\left[\left(\mathbb{A}(\geq
1|_p)+\sum_{k=1}^ni^ka_k\right)\bmod p\right]$.

In addition, the first necessary and sufficient condition ensures
that the range of $F_{0,0},\cdots,F_{0,p-1}$ forms a partition of
$\mathbb{A}\backslash\mathbb{A}(\geq 1|_p)$. This means that there
exists a super-bijection $F_0$ over
$\mathbb{A}\backslash\mathbb{A}(\geq 1|_p)$.

The above analyses show that $f(x)$ forms a bijection over
$\mathbb{A}$.

2) Assume the ``if" part is true for $2,\cdots,d-1$. Let us prove
the case of $d\geq 3$. Similarly, let us consider the $p$
sub-bijections, $F_1$, $F_{0,1},\cdots,F_{0,p-1}$, separately.

When $x\in\mathbb{A}(\geq 1|_p)$, assume $x=pz$, where
$z\in\{0,\cdots,p^{d-1}-1\}$, so
$f(x)=a_nx^n+\cdots+a_2x^2+a_1x\equiv
p(a_np^{n-1}x^n+\cdots+a_2pz+a_1z)\pmod{p^d}$. Apparently, $f(x)$ is
uniquely determined by the polynomial
$f_1^*(z)=b_nz^n+\cdots+b_1z=a_np^{n-1}x^n+\cdots+a_2pz+a_1z$ modulo
$p^{d-1}$. Since $b_i\equiv 0\pmod p$ when $i\geq 2$, one can easily
verify that $f_1^*(z)$ satisfies the two necessary and sufficient
conditions, so from the previous assumption, $f_1^*(z)$ is a
permutation polynomial modulo $p^{d-1}$. This means that $f(x)$
forms a bijection over $\mathbb{A}(\geq 1|_p)$.

$\forall i=1\sim p-1$ and $\forall x\in[(\mathbb{A}(\geq
1|_p)+i)\bmod p]$, assume $x=y+i$, where $y\in\mathbb{A}(\geq
1|_p)$. Substitute $x=y+i$ into $f(x)$, one has
$f_{0,i}^*(y)=f_{0,i}^{**}(y)+\sum_{k=1}^ni^ka_k$, where
$f_{0,i}^{**}(y)=\sum_{l=n}^1\left(\sum_{k=l}^n\binom{k}{l}i^{k-l}a_ky^l\right)$.
Assume $y=pz$, where $z\in\{0,\cdots,p-1\}$, one has
$f_{0,i}^{**}(y)=\sum_{l=n}^1\left(\sum_{k=l}^n\binom{k}{l}i^{k-l}a_kp^lz^l\right)
=p\sum_{l=n}^1\left(\sum_{k=l}^n\binom{k}{l}i^{k-l}a_kp^{l-1}z^l\right)$.
Apparently, $f_{0,i}^{**}(y)$ is uniquely determined by the
polynomial
$f_{0,i}^{***}(z)=b_nz^n+\cdots+b_1z=\sum_{l=n}^2\left(\sum_{k=l}^n\binom{k}{l}i^{k-l}a_kp^{l-1}z^l\right)+
\sum_{k=1}^nki^{k-1}a_kz$ modulo $p^{d-1}$. Similarly, since
$b_k\equiv 0\pmod p$ when $k\geq 2$, $\forall
i,j\in\{0,\cdots,p-1\}$ and $i\neq j$, one has
$\sum_{k=1}^nki^{k-1}b_k\equiv b_1=\sum_{k=1}^nki^{k-1}a_k\not\equiv
0\pmod p$ and $\sum_{k=1}^nb_k(j^k-i^k)\equiv
b_1(j-i)=(j-i)\sum_{k=1}^nki^{k-1}a_k\not\equiv 0\pmod p$, where
note that $j-i\not\equiv 0\pmod p$. That is, the two necessary and
sufficient conditions hold for $f_{0,i}^{***}(z)$, so from the
previous assumption, $f_{0,i}^{***}(z)$ is a permutation polynomial
modulo $p^{d-1}$, i.e., $f_{0,i}^{**}(y)$ forms a bijection over
$\mathbb{A}(\geq 1|_p)$ and $f(x)$ forms a bijection
$F_{0,i}:[(\mathbb{A}(\geq 1|_p)+i)\bmod
p]\to\left[\left(\mathbb{A}(\geq
1|_p)+\sum_{k=1}^ni^ka_k\right)\bmod p\right]$.

In addition, the first necessary and sufficient condition ensures
that the range of $F_{0,0},\cdots,F_{0,p-1}$ forms a partition of
$\mathbb{A}\backslash\mathbb{A}(\geq 1|_p)$. This means that there
exists a super-bijection $F_0$ over
$\mathbb{A}\backslash\mathbb{A}(\geq 1|_p)$.

The above analyses show that $f(x)$ forms a bijection over
$\mathbb{A}$. Thus this theorem is proved.
\end{proof}

\begin{corollary}\label{corollary:bijective-poly-NSC-2d}
The polynomial $f(x)=a_nx^n+\cdots+a_1x$ is a permutation polynomial
modulo $2^d$ if and only if the following two conditions are true
simultaneously: $a_1\equiv 1\pmod 2$, $a_2+a_4+\cdots\equiv
a_3+a_5+\cdots\equiv 0\pmod 2$.
\end{corollary}
\begin{proof}
From Theorem \ref{theorem:bijective-poly-NSC-pd}, choosing $p=2$,
one has the following necessary and sufficient conditions:
$\sum_{i=1}^na_i\not\equiv 0\pmod 2$, $a_1\not\equiv 0\pmod 2$ and
$\sum_{i=1}^nia_i\not\equiv 0\pmod 2$. These conditions can be
simplified to be: $a_1\equiv 1\pmod 2$, $\sum_{i=2}^na_i\equiv
\sum_{i=2}^nia_i\equiv 0\pmod 2$. Removing even terms from
$\sum_{i=2}^nia_i\equiv 0\pmod 2$, one has $a_3+a_5+\cdots\equiv
0\pmod 2$. Then, subtracting $a_3+a_5+\cdots$ from
$\sum_{i=2}^na_i$, one has $a_2+a_4+\cdots\equiv 0\pmod 2$. This
corollary is thus proved.
\end{proof}

\begin{corollary}\label{corollary:bijective-poly-NSC-all-degrees}
Assume $p$ is a prime and $d\geq 2$. If $f(x)$ is a permutation
polynomial modulo $p^d$, then $\forall i\geq 1$, it is still a
permutation polynomial modulo $p^i$.
\end{corollary}
\begin{proof}
This corollary is a direct result of Theorem
\ref{theorem:bijective-poly-NSC-pd}.
\end{proof}

\begin{theorem}\label{theorem:bijective-poly-NSC-pd-degree2}
Assume $p$ is a prime and $d\geq 1$. A binomial $f(x)=a_2x^2+a_1x$
is a permutation polynomial modulo $p^d$ if and only if
$a_1\not\equiv 0\pmod p$ and $a_2\equiv 0\pmod p$.
\end{theorem}
\begin{proof}
The ``if" part can be easily verified by checking the necessary and
sufficient conditions in Theorem
\ref{theorem:bijective-poly-NSC-pd}. We focus on the ``only if"
part.

When $p=2$, one can verify the result is true. Let us consider the
case of $p>2$. From Theorem \ref{theorem:bijective-poly-NSC-pd},
$a_1(1-(p-1))+a_2(1^2-(p-1)^2)\equiv 2a_1\not\equiv 0\pmod p$, which
immediately leads to $a_1\not\equiv 0\pmod p$. Again, from Theorem
\ref{theorem:bijective-poly-NSC-pd}, the following conditions hold:
$\forall i=1\sim p-1$, $ia_1+i^2a_2\not\equiv 0\pmod p$. Since $p$
is a prime, each integer in $\{1,\cdots,p-1\}$ has an inverse modulo
$p$. Multiplying the inverse of $i$ at both sides of each condition,
one gets $a_1+ia_2\not\equiv 0\pmod p$, so $a_2\not\equiv
-a_1\bar{i}\pmod p$, where $\bar{i}$ is the inverse of $i$ modulo
$p$. Since $\{i\}$ forms a reduced system of residues modulo $p$,
$\{-a_1\bar{i}\}=\{-a_1,\cdots,-a_1(p-1)\}$ still forms a reduced
system of residues modulo $p$. Thus $a_2\equiv 0\pmod p$.
\end{proof}
\begin{remark}
Note that Theorem \ref{theorem:bijective-poly-NSC-pd-degree2}
actually says that the first group of the necessary and sufficient
conditions covers the second group of conditions when the degree is
1 or 2 modulo $p^d$.
\end{remark}
\begin{corollary}
Assume $p$ is a prime and $d\geq 1$. If $f(x)$ is a permutation
polynomial of degree 1 modulo $p^d$, then $\forall i\geq 1$, it is
still a permutation polynomial modulo $p^i$.
\end{corollary}
\begin{proof}
This corollary is a direct result of Theorem
\ref{theorem:bijective-poly-NSC-pd-degree2}.
\end{proof}

\fproblem{When $3\leq n\leq p-1$, is it possible to get further
simplified necessary and sufficient conditions?}

\begin{remark}\label{remark:RepeatingWork}
After finishing the first draft of this paper, we noticed Rivest's
paper \cite{Rivest:PPmod2w:FFTA2001} and found Corollary
\ref{corollary:bijective-poly-NSC-2d} was proved by the author in
2002. Through \cite{Rivest:PPmod2w:FFTA2001}, we further noticed
Mullen's paper \cite{Mullen:PolyFun-mod:AMH1984} and realized that
Theorem \ref{theorem:bijective-poly-NSC-pd} can also be derived from
Theorem 123 in \cite{HW:NumberTheory1979}, where the second
condition becomes that $f'(x)\not\equiv 0\pmod p$ holds for any
integer $x$. In addition, recently we found yet another paper
\cite{Sun:PPmodCoding:IEEETIT2005}, in which Theorem
\ref{theorem:bijective-poly-NSC-pd-degree2} was also obtained in a
similar way (Corollary 2.4). Furthermore, we also noticed Lemma 4.2
in Chap. 4 of \cite{RMT:DicksonPoly1993} gives a more general form
of Theorem \ref{theorem:bijective-poly-NSC-pd}. Considering the fact
that our proof of Theorem \ref{theorem:bijective-poly-NSC-pd} is
independent of Theorem 123 in \cite{HW:NumberTheory1979}, it can be
considered as a different proof of this result.
\end{remark}

\subsection{Counting Permutation Polynomials and Induced Permutations}

The case of $n\geq p-1$ modulo $p$ has been solved in Corollary
\ref{corollary:count-bijective-poly-p}. This subsection discusses
other cases modulo $p^d$ ($d\geq 1$).

\begin{notation}
Assume $p$ is a prime and $d\geq 1$. Denote the number of
permutation polynomials and the number of all polynomials in a
complete system of polynomial resides of degree $\leq n$ modulo
$p^d$ by $N_{pp}(\leq n,p^d)$ and $N_p(\leq n,p^d)$, respectively.
\end{notation}

\begin{remark}
Note that the number of permutation polynomials of degree $n$ modulo
$p^d$ can be easily calculated to be $N_{pp}(\leq n,p^d)-N_{pp}(\leq
n-1,p^d)$. So this paper only focuses on the number of permutation
polynomials of degree $\leq n$ modulo $p^d$.
\end{remark}

\begin{theorem}
For any prime $p$ and $d\geq 1$, $\dfrac{N_{pp}(\leq
1,p^d)}{N_p(\leq 1,p^d)}=\dfrac{p-1}{p}$ and $\dfrac{N_{pp}(\leq
2,p^d)}{N_p(\leq 2,p^d)}=\dfrac{p-1}{p^2}$.
\end{theorem}
\begin{proof}
This theorem is a direct result of Theorems
\ref{theorem:bijection-degree1} and
\ref{theorem:bijective-poly-NSC-pd-degree2}.
\end{proof}

\begin{theorem}\label{theorem:count-bijective-poly-pd}
For any prime $p$ and $d\geq 2$, $\dfrac{N_{pp}(\leq
n,p^d)}{N_p(\leq n,p^d)}=\dfrac{(p-1)^p(p-1)!}{p^{2p-1}}$ when
$n\geq 2p-1$.
\end{theorem}
\begin{proof}
From Theorem \ref{theorem:bijective-poly-NSC-pd}, a bijective
polynomial should satisfy the following conditions:
\begin{itemize}
\item
$\binom{p}{2}=\frac{p(p-1)}{2}$ conditions: $\forall
i,j\in\{0,\cdots,p-1\}$ and $i\neq j$,
$\sum_{k=1}^na_k(i^k-j^k)=a_1(i-j)+\cdots+a_n(i^n-j^n)\not\equiv
0\pmod p$;

\item
$p$ conditions: $\forall i\in\{0,\cdots,p-1\}$,
$\sum_{k=1}^nki^{k-1}a_k=a_1+2ia_2+\cdots+ni^{n-1}a_n\not\equiv
0\pmod p$.
\end{itemize}
Among the above $\binom{p}{2}+p=\frac{p(p+1)}{2}$ conditions, choose
the following $2p-1$ conditions:
\begin{itemize}
\item
$p-1$ conditions: $\forall i\in\{1,\cdots,p-1\}$ and $j=0$,
$\sum_{k=1}^na_ki^k=a_1i+\cdots+a_ni^n\equiv b_i\pmod p$, where
$b_i\not\equiv 0\pmod p$;

\item
$p$ conditions: $\forall i\in\{0,\cdots,p-1\}$,
$\sum_{k=1}^nki^{k-1}a_k=a_1+2ia_2+\cdots+ni^{n-1}a_n\equiv
b_{p+i}\pmod p$, where $b_{p+i}\not\equiv 0\pmod p$.
\end{itemize}
Rewrite the above $2p-1$ condition as a system of congruences:
\[
\left[\begin{matrix}%
1 & 1 & 1 & \cdots & 1\\
2 & 2^2 & 2^3 & \cdots & 2^n\\
\vdots & \vdots & \vdots & \ddots & \vdots\\
(p-1) & (p-1)^2 & (p-1)^3 & \cdots & (p-1)^n\\
1 & 0 & 0 & \cdots & 0\\
1 & 2 & 3 & \cdots & n\\
1 & 2\cdot 2 & 3\cdot 2^2 & \cdots & n\cdot 2^{n-1}\\
\vdots & \vdots & \vdots & \ddots & \vdots\\
1 & 2\cdot(p-1) & 3\cdot(p-1)^2 & \cdots & n\cdot(p-1)^{n-1}
\end{matrix}\right]
\left[\begin{matrix}%
a_1\\a_2\\a_3\\a_4\\\vdots\\a_{n-3}\\a_{n-2}\\a_{n-1}\\a_n
\end{matrix}\right]\equiv
\left[\begin{matrix}%
b_1\\
b_2\\
\vdots\\
b_{p-1}\\
b_p\\
b_{p+1}\\
b_{p+2}\\
\vdots\\
b_{2p-1}
\end{matrix}\right]\pmod p.
\]
If we only consider $a_1,\cdots,a_{2p-1}$ as unknown variables, the
above system can be reduced to be the following system:\small
\[
\left[\begin{matrix}%
1 & 1 & 1 & \cdots & 1\\
2 & 2^2 & 2^3 & \cdots & 2^{2p-1}\\
\vdots & \vdots & \vdots & \ddots & \vdots\\
(p-1) & (p-1)^2 & (p-1)^3 & \cdots & (p-1)^{2p-1}\\
1 & 0 & 0 & \cdots & 0\\
1 & 2 & 3 & \cdots & 2p-1\\
1 & 2\cdot 2 & 3\cdot 2^2 & \cdots & (2p-1)\cdot 2^{2p-2}\\
\vdots & \vdots & \vdots & \ddots & \vdots\\
1 & 2\cdot(p-1) & 3\cdot(p-1)^2 & \cdots & (2p-1)\cdot(p-1)^{2p-2}
\end{matrix}\right]
\left[\begin{matrix}%
a_1\\a_2\\\vdots\\a_{p-1}\\a_p\\a_{p+1}\\a_{p+2}\\\vdots\\a_{2p-1}
\end{matrix}\right]\equiv
\left[\begin{matrix}%
b_1-\sum_{k=2p}^na_i\\
b_2-\sum_{k=2p}^na_i2^k\\
\vdots\\
b_{p-1}-\sum_{k=2p}^na_i(p-1)^k\\
b_p\\
b_{p+1}-\sum_{k=2p}^nka_i\\
b_{p+2}-\sum_{k=2p}^nk2^{k-1}a_i\\
\vdots\\
b_{2p-1}-\sum_{k=2p}^nk(p-1)^{k-1}a_i
\end{matrix}\right]\pmod p.
\]\normalsize
Denoting the above system by $\bm{A}\bm{X}\equiv\bm{B}\pmod p$, from
Corollary \ref{corollary:det-powers}, one has
\[
|\bm{A}|=(-1)^{\frac{(p-1)(p-2)}{2}}\prod_{i=1}^{p-1}i^4\prod_{1\leq
i<j\leq p-1}(j-i)^4.
\]
Since all factors of $|\bm{A}|$ are in $\{1,\cdots,p-1\}$ and $p$ is
a prime, $\gcd(|\bm{A}|,p)=1$. Then, for each valid combination of
$(b_1,\cdots,b_{2p-1},a_{2p},\cdots,a_n)$, the above system of
congruences has a unique set of incongruent solutions.

Next, let us count the number of all valid combinations of
$(b_1,\cdots,b_{2p-1},a_{2p},\cdots,a_n)$. It is obvious that
$\{a_{2p},\cdots,a_n\}$ can be any value and $b_p,\cdots,b_{2p-1}$
can be any nonzero value modulo $p$. However, $b_1,\cdots,b_{p-1}$
are also constrained by the following conditions: $\forall
i,j\in\{1,\cdots,p-1\}$ and $i\neq j$,
$\sum_{k=1}^na_k(i^k-j^k)\equiv b_i-b_j\not\equiv 0\pmod p$. That
is, $\{b_i\bmod p\}_{i=1}^{p-1}$ forms a complete permutation over
$\{1,\cdots,p-1\}$, so the number of possible values of
$(b_1,\cdots,b_{p-1})$ is $(p-1)!$ in total $p^{p-1}$ combinations
of the $p-1$ values. Combining the above fact, one immediately gets
$N_{pp}(\leq n,p^d)/N_p(\leq n,p^d)=(1-1/p)^p\cdot
(p-1)!/p^{p-1}=(p-1)^p(p-1)!/p^{2p-1}$.

Thus this theorem is proved.
\end{proof}

\begin{theorem}
For any prime $p$ and $d\geq 2$, the following inequalities hold:
\begin{enumerate}
\item
when $3\leq n\leq p$, $\dfrac{N_{pp}(\leq n,p^d)}{N_p(\leq
n,p^d)}\leq\dfrac{(p-1)P(p-1,n-1)}{p^n}=\dfrac{(p-1)(n-1)!\binom{p-1}{n-1}}{p^n}$;

\item
when $p+1\leq n\leq 2p-2$, $\dfrac{N_{pp}(\leq n,p^d)}{N_p(\leq
n,p^d)}\leq\dfrac{(p-1)!}{p^{p-1}}\left(\dfrac{p-1}{p}\right)^n$.
\end{enumerate}
\end{theorem}
\begin{proof}
When $3\leq n\geq 2p-2$, the matrix in the proof of Theorem
\ref{theorem:count-bijective-poly-pd} has at most $n$ free
congruences and other $2p-1-n$ congruences are actually linear
combinations of the $n$ free ones. This means that there exists an
upper bound of $\frac{N_{pp}(\leq n,p^d)}{N_p(\leq n,p^d)}$. Note
that the values of $b_1,\cdots,b_{p-1}$ should form a permutation
over $\{1,\cdots,p-1\}$, so we consider the following two
conditions, respectively.

1) When $3\leq n\leq p$, $b_p$ has $(p-1)$ possible values and the
$n-1$ left free variables have $P(p-1,n-1)=(n-1)!\binom{p-1}{n-1}$
combinations, so $\dfrac{N_{pp}(\leq n,p^d)}{N_p(\leq
n,p^d)}\leq\dfrac{(p-1)P(p-1,n-1)}{p^n}=\dfrac{(p-1)(n-1)!\binom{p-1}{n-1}}{p^n}$.

2) When $p+1\leq n\leq 2p-2$, $p-1$ free variables form the
permutation over $\{1,\cdots,p-1\}$ and other $n-(p-1)$ variables
are totally free, so the number of possibilities of the $n$ free
variables is $(p-1)!(p-1)^{n-(p-1)}$. Thus, $\dfrac{N_{pp}(\leq
n,p^d)}{N_p(\leq
n,p^d)}\leq\dfrac{(p-1)!(p-1)^{n-(p-1)}}{p^n}=\dfrac{(p-1)!}{p^{p-1}}\left(\dfrac{p-1}{p}\right)^n$.
\end{proof}

\begin{theorem}
For any prime $p$ and $3\leq n\leq p-2$, the following inequality
hold:
\[
\dfrac{N_{pp}(\leq n,p)}{N_p(\leq
n,p)}\leq\dfrac{P(p-1,n-1)}{p^n}=\dfrac{(n-1)!\binom{p-1}{n-1}}{p^n}.
\]
\end{theorem}
\begin{proof}
When $d=1$, the second group of necessary and sufficient conditions
disappear. Then, following the similar idea of proving the above
theorem, this theorem is proved.
\end{proof}

\fproblem{When $3\leq n\leq 2p-2$, it is still possible to get a
close form of the \textbf{exact} value of $\dfrac{N_{pp}(\leq
n,p^d)}{N_p(\leq n,p^d)}$?}

\begin{example}
When $d\geq 2$,
\[
\frac{N_{pp}(\leq n,2^d)}{N_p(\leq n,2^d)}=1/2^{\min(n,3)}=\begin{cases}%
1/2, & n=1,\\
1/2^2, & n=2,\\
1/2^3, & n\geq 3.
\end{cases}
\]
\end{example}
\begin{solution}
When $n=1$, from Lemma \ref{theorem:bijection-degree1}, a polynomial
is a permutation polynomial modulo $2^d$ if and only if
$\gcd(a_1,2^d)=1$. So, $a_1\equiv 1\pmod 2$, which means
$\frac{N_{pp}(\leq 1,2^d)}{N_p(\leq 1,2^d)}=1/2$.

When $n=2$, assume $f(x)=a_2x^2+a_1x$. From Theorem
\ref{theorem:bijective-poly-NSC-pd-degree2}, the necessary and
sufficient conditions are $a_1\equiv 1\pmod 2$ and $a_2\equiv 0\pmod
2$. Thus, $\frac{N_{pp}(\leq 2,2^d)}{N_p(\leq 2,2^d)}=1/2^2$.

When $n\geq 2p-1=3$, from Theorem
\ref{theorem:count-bijective-poly-pd}, one has $\frac{N_{pp}(\leq
n,2^d)}{N_p(\leq n,2^d)}=(2-1)^3\cdot(2-1)!/2^{2\cdot 2-1}=1/2^3$.

Computer experiments have been made to verify the above results.
\end{solution}

\begin{example}
When $d\geq 2$,
\[
\frac{N_{pp}(\leq n,3^d)}{N_p(\leq n,3^d)}=\begin{cases}%
2/3, & n=1,\\
2/3^2, & n=2,\\
4/3^n, & n=3,4,\\
16/3^5, & n\geq 5.
\end{cases}
\]
\end{example}
\begin{solution}
When $n=1$, from Lemma \ref{theorem:bijection-degree1}, a polynomial
is a permutation polynomial modulo $3^d$ if and only if
$\gcd(a_1,3^d)=1$. So, $a_1\equiv 1,2\pmod 3$, which means
$\frac{N_{pp}(\leq 1,3^d)}{N_p(\leq 1,3^d)}=2/3$.

When $n=2$, assume $f(x)=a_2x^2+a_1x$. From Theorem
\ref{theorem:bijective-poly-NSC-pd-degree2}, the necessary and
sufficient conditions are $a_1\not\equiv 0\pmod 3$ and $a_2\equiv
0\pmod 3$. This means that $\frac{N_{pp}(\leq 2,3^d)}{N_p(\leq
2,3^d)}=2/3^2$.

When $n=3$, from Theorem \ref{theorem:bijective-poly-NSC-pd}, the
necessary and sufficient conditions are as follows:
$a_1+a_2+a_3\not\equiv 0\pmod 3$, $2a_1+4a_2+8a_3\not\equiv 0\pmod
3$, $(2-1)a_1+(4-1)a_2+(8-1)a_3\not\equiv 0\pmod 3$, $a_1\not\equiv
0\pmod 3$, $a_1+2a_2+3a_3\not\equiv 0\pmod 3$ and
$a_1+4a_2+12a_3\not\equiv 0\pmod 3$. These conditions can be further
simplified as $a_1\not\equiv 0\pmod 3$, $a_2\equiv 0\pmod 3$ and
$a_1+a_3\not\equiv 0\pmod 3$. So, the possible values of
$(a_1,a_2,a_3)$ modulo 3 are $(1,0,0)$, $(1,0,1)$, $(2,0,0)$ and
$(2,0,2)$. Thus, $\frac{N_{pp}(\leq 3,3^d)}{N_p(\leq 3,3^d)}=4/3^3$.
In the same way, one can deduce the results when $n=4$.

When $n\geq 2p-1=5$, from Theorem
\ref{theorem:count-bijective-poly-pd}, one has $\frac{N_{pp}(\leq
n,3^d)}{N_p(\leq n,3^d)}=(3-1)^3\cdot(3-1)!/3^{2\cdot 3-1}=16/3^5$.

Computer experiments have been made to verify the above results.
\end{solution}

After getting the number of permutation polynomials of degree $\leq
n$ modulo $p^d$, one can easily calculate the number of distinct
permutations induced by the permutation polynomials of degree $\leq
n$ modulo $p^d$, by using Lemma \ref{lemma:Equ-Poly-Null-Poly} and
the results on null polynomials modulo $p^d$ given in
\cite{Li:NullPoly2005}. We have the following theorem.

\begin{theorem}
Assume $p$ is a prime, $d\geq 1$ and $N_{np}(\leq n,p^d)$ denotes
the number of null polynomials of degree $\leq n$ modulo $p^d$.
Then, the number of distinct permutations induced from polynomials
of degree $\leq n$ modulo $p^d$ is $N_{pp}(\leq n,p^d)/N_{np}(\leq
n,p^d)$.
\end{theorem}
\begin{proof}
It is obvious since each polynomial has $N_{np}(\leq n,p^d)$
equivalent polynomials from Lemma \ref{lemma:Equ-Poly-Null-Poly}.
\end{proof}

\begin{remark}
Note that Corollary 4.1 of \cite{KO:CountingPF:DMJ1968} gives a
different proof of Theorem \ref{theorem:count-bijective-poly-pd}.
However, \cite{KO:CountingPF:DMJ1968} mainly focuses on the total
number of non-equivalent polynomial functions of arbitrary degree
modulo $p^d$ and does not study the case when the degree is also
given. So, the results given in this paper are more complete.
\end{remark}

\subsection{Determining (Permutation) Polynomials from Induced Bijection}

In this subsection, we study the problem of determining all
equivalent polynomials when the induced polynomial function is (or
partially) known. Note that the following results are also valid for
polynomials that induce any polynomial functions (maybe not
bijections).

Following Lemma \ref{lemma:Equ-Poly-Null-Poly} and the results
obtained in \cite{Li:NullPoly2005}, once we get one permutation
polynomials inducing the given bijection, we can determine all
equivalent permutation polynomials. So, it is sufficient to derive
only one equivalent polynomial as a seed.

\begin{theorem}\label{theorem:solving-poly-small-n}
Assume $p$ is a prime, $d\geq 1$ and $f(x)=a_nx^n+\cdots+a_1x+a_0$
is a polynomial of degree $n\leq p-1$ modulo $p^d$. Given
$x_0,\cdots,x_n\in\mathbb{Z}$, if $\forall i,j\in\{0,\cdots,n\}$ and
$i\neq j$, $x_i\not\equiv x_j\pmod p$, then $f(x)$ can be uniquely
determined by solving the following system of congruence:
\begin{equation}
\left[\begin{matrix}%
1 & x_0 & x_0^2 & \cdots & x_0^n\\
1 & x_1 & x_1^2 & \cdots & x_1^n\\
\vdots & \vdots & \vdots & \ddots & \vdots\\
1 & x_n & x_n^2 & \cdots & x_n^n
\end{matrix}\right]%
\left[\begin{matrix}%
a_0\\a_1\\\vdots\\a_n
\end{matrix}\right]\equiv
\left[\begin{matrix}%
f(x_0)\\f(x_1)\\\vdots\\f(x_n)
\end{matrix}\right]\pmod{p^d}.\label{equation:solve-PP-nLessp}
\end{equation}
\end{theorem}
\begin{proof}
Denote the system of congruences by
$\bm{A}_x\bm{X}_a\equiv\bm{B}_f\pmod{p^d}$. Since $\bm{A}_x$ is a
Vondermonde matrix, $|\bm{A}_x|=\prod_{0\leq i<j\leq n}(x_j-x_i)$.
From $x_i\not\equiv x_j\pmod p$, $\gcd(x_j-x_i,p)=1$, so
$\gcd(|\bm{A}_x|,p^d)=1$. Thus, the system of congruences has a
unique set of incongruent solutions and this theorem is proved.
\end{proof}
\begin{remark}
When $n=p-1$, it is obvious that $x_0,\cdots,x_n$ form a complete
system of residues modulo $p$. When $n<p-1$, $x_0,\cdots,x_n$ form
an incomplete system of residues modulo $p$. The simplest choice of
the $n+1$ values is: $\{x_i=i\}_{i=0}^n=\{0,\cdots,n\}$.
\end{remark}
\begin{corollary}\label{corollary:solving-poly-small-n}
Assume $p$ is a prime, $d\geq 1$ and $f(x)=a_nx^n+\cdots+a_1x$ is a
polynomial of degree $n\leq p-1$ modulo $p^d$. Given
$x_1,\cdots,x_n\not\equiv 0\pmod p$, if $\forall
i,j\in\{1,\cdots,n\}$ and $i\neq j$, $x_i\not\equiv x_j\pmod p$,
then $f(x)$ can be uniquely determined by solving the following
system of congruence:
\begin{equation}
\left[\begin{matrix}%
x_1 & x_1^2 & \cdots & x_1^n\\
x_2 & x_2^2 & \cdots & x_2^n\\
\vdots & \vdots & \ddots & \vdots\\
x_n & x_n^2 & \cdots & x_n^n
\end{matrix}\right]%
\left[\begin{matrix}%
a_1\\a_2\\\vdots\\a_n
\end{matrix}\right]\equiv
\left[\begin{matrix}%
f(x_1)\\f(x_2)\\\vdots\\f(x_n)
\end{matrix}\right]\pmod{p^d}.
\end{equation}
\end{corollary}
\begin{proof}
This corollary is a special case of Theorem
\ref{theorem:solving-poly-small-n}. For $i=1\sim n$, factoring out
$x_i$ from row $i$ of the matrix, one immediately has
$|\bm{A}|=\prod_{i=1}^nx_i\prod_{1\leq i<j\leq n}(x_j-x_i)$. From
the conditions of $\{x_i\}_{i=1}^n$, $|\bm{A}|$ is relatively prime
to $p$ and the system of congruences has a unique set of incongruent
solutions, thus this corollary is proved.
\end{proof}

When $n\geq p$ or the value of $n$ is unknown, the above method
cannot be directly used to determine polynomials that induce the
given polynomial function. If we can find a way to reduce the degree
of polynomials, then the above method can be employed to determine
the coefficients. In the following, we give a way to achieve this
task.

\begin{lemma}\label{lemma:determine-PP-dLessp}
Assume $p$ is a prime and $2\leq d\leq p$. If
$f(x)=a_nx^n+\cdots+a_1x+a_0$ is a polynomial of degree $n\leq pd-1$
modulo $p^d$, then all its equivalent polynomials $\leq pd-1$ modulo
$p^d$ can be determined from the induced polynomial function over
$\{0,\cdots,p^d-1\}$.
\end{lemma}
\begin{proof}
Choosing $x=py_1+b_0$, i.e., $y_1=\lfloor x/y\rfloor$ and
$b_0=(x\bmod p)\in\{0,\cdots,p-1\}$, we have $p$ sub-polynomials:
\[
f_{b_0}(y_1)=\begin{cases}%
\sum_{k=n}^0p^ka_ky_1^k, & \mbox{when }b_0=0,\\
\sum_{k=n}^0\left(\sum_{l=n}^ka_l\binom{l}{k}p^ky_1^kb_0^{l-k}\right),
& \mbox{when }b_0\in\{1,\cdots,p-1\}.
\end{cases}
\]
Apparently, they have a uniform form:
$f_{b_0}(y_1)=\sum_{k=n}^0p^ka_k^{(b_0)}y_1^k$, where
$a_k^{(0)}=a_k$ and
$a_k^{(b_0)}=\sum_{l=n}^ka_l\binom{l}{k}b_0^{l-k}$ when
$b_0\in\{1,\cdots,p-1\}$. Note that $f_{b_0}(y_1)\equiv
f_{b_0}^*(y_1)=\sum_{k=d-1}^0p^ka_k^{(b_0)}y_1^k\pmod{p^d}$, so we
can focus on $f_{b_0}^*(y_1)$ only. Since $d\leq p$, each
sub-polynomial $f_{b_0}^*(y_1)$ is of degree less than $p$ modulo
$p^d$, so all the coefficients can be uniquely solved modulo $p^d$,
i.e. the value of each $a_k^{(b_0)}$, can be uniquely solved modulo
$p^{d-k}$. It is obvious that the $p^k$ distinct valid values of
$a_k^{(b_0)}$ modulo $p^d$ are equivalent for $f_{b_0}(y_1)$ and so
equivalent for $f(x)$.

For each set of the valid values of all coefficients modulo $p^d$,
one has a system of congruences in the form
$\bm{A}\bm{X}_a\equiv\bm{B}\pmod{p^d}$:
\begin{equation}
\left[\begin{matrix}%
\bm{A}_0\\\bm{A}_1\\\vdots\\\bm{A}_{p-1}
\end{matrix}\right]
\left[\begin{matrix}%
a_0\\a_1\\\vdots\\a_{pd-1}
\end{matrix}\right]\equiv
\left[\begin{matrix}%
\bm{B}_0\\\bm{B}_1\\\vdots\\\bm{B}_{p-1}
\end{matrix}\right]\pmod{p^d},\label{equation:solve-PP-pd}
\end{equation}
where
\[
\bm{A}_0=\left[\begin{matrix}%
\bm{I}_{d\times d} & \bm{0}_{d \times d(p-1)}\end{matrix}\right]
=\left[\begin{array}{ccccc:ccccc}%
1 & 0 & 0 & \cdots & 0 & 0 & 0 & 0 & \cdots & 0\\
0 & 1 & 0 & \cdots & 0 & 0 & 0 & 0 & \cdots & 0\\
0 & 0 & 1 & \cdots & 0 & 0 & 0 & 0 & \cdots & 0\\
\vdots & \vdots & \vdots & \ddots & \vdots & \vdots & \vdots &
\vdots & \ddots & \vdots\\
0 & 0 & 0 & \cdots & 1 & 0 & 0 & 0 & \cdots & 0
\end{array}\right]_{d\times pd},
\]
for $b_0=1\sim p-1$,
\[
\bm{A}_{b_0}=\left[\begin{matrix}%
\bm{A}_{b_0}^{(L)} & \bm{A}_{b_0}^{(R)}\end{matrix}\right]=
\left[\begin{array}{ccccc:ccc}%
1 & b_0 & b_0^2 & \cdots & b_0^{d-1} & b_0^d & \cdots & b_0^{pd}\\
0 & 1 & 2b_0 & \cdots & (d-1)b_0^{d-2} & db_0^{d-1} & \cdots & (pd-1)b_0^{pd-1}\\
0 & 0 & \binom{2}{2} & \cdots & \binom{d-1}{2}b_0^{d-3} &
\binom{d}{2}b_0^{d-2} & \cdots & \binom{pd-1}{2}b_0^{pd-2}\\
\vdots & \vdots & \vdots & \ddots & \vdots & \vdots & \ddots & \vdots\\
0 & 0 & 0 & \cdots & \binom{d-1}{d-1} & \binom{d}{d-1}b_0 & \cdots &
\binom{pd-1}{d-1}b_0^{pd-(d-1)}
\end{array}\right]_{d\times pd},
\]
and for $b_0=0\sim p-1$,
$\bm{B}_{b_0}=\left[\begin{matrix}a_0^{(b_0)} & a_1^{(b_0)} &
a_{d-1}^{(b_0)}\end{matrix}\right]^T$. From Lemma
\ref{lemma:det-binom-powers}, one can see $|\bm{A}|$ is relatively
prime to $p$. Thus, for each valid set of the values of the
coefficients $\left\{a_k^{(b_0)}\right\}_{0\leq k\leq p-1 \atop
0\leq b_0\leq p-1}$, the above system of congruences has a unique
set of solutions modulo $p^d$. One can easily verify that each set
of solutions corresponds to an equivalent of the polynomial $f(x)$.
Thus this lemma is proved.
\end{proof}
\begin{remark}
In fact, in the proof of the above lemma, we can also calculate the
number of equivalent polynomials of $f(x)$ of degree $\leq pd-1$
modulo $p^d$. It is $p^{(1+\cdots+(d-1))p}=p^{\frac{d(d-1)p}{2}}$.
From Lemma \ref{lemma:Equ-Poly-Null-Poly}, this number should be
equal to the number of null polynomials of degree $\leq pd-1$ modulo
$p^d$. Clearly, $p^{\frac{d(d-1)p}{2}}$ agrees with the results
(Lemma 34 and Theorem 43) obtained in \cite{Li:NullPoly2005}.
\end{remark}

\begin{corollary}
Assume $p$ is a prime and $2\leq d\leq p$. If
$f(x)=a_nx^n+\cdots+a_1x+a_0$ is a polynomial modulo $p^d$, then all
its equivalent polynomials modulo $p^d$ an be determined from the
induced polynomial function over $\{0,\cdots,p^d-1\}$.
\end{corollary}
\begin{proof}
When $n\geq pd$, one can move $a_{pd},\cdots,a_n$ to the right side
of the matrix in the proof of the above lemma. Then,
$a_{pd},\cdots,a_n$ become free variables, so each set of their
values corresponds to $p^{\frac{d(d-1)p}{2}}$ equivalent polynomials
modulo $p^d$. That is, in total we have
$p^{d(n-pd+1)\frac{d(d-1)p}{2}}$ equivalent polynomials. Thus this
corollary is also true.
\end{proof}

\begin{theorem}
Assume $p$ is a prime and $d\geq 2$. If
$f(x)=a_nx^n+\cdots+a_1x+a_0$ is a polynomial modulo $p^d$, then all
its equivalent polynomials modulo $p^d$ an be determined from the
induced polynomial function over $\{0,\cdots,p^d-1\}$.
\end{theorem}
\begin{proof}
We use induction on $d$ to prove this theorem. The case of $2\leq
d\leq p$ has been proved above. Let us prove the case of $d>p$ under
the assumption that this theorem is true for any integer less than
$d$.

Using the same way in the proof of Lemma
\ref{lemma:determine-PP-dLessp}, we can get $p$ sub-polynomials
$f_{b_0}^*(y_1)$, which uniquely determine the induced bijection.
Since the degree of $f_{b_0}^*(y_1)$ modulo $p^d$ is not less than
$p$, the coefficients cannot be uniquely solved, let us try to
further decompose each sub-polynomial in the same way.

At first, note that the value of $a_0^{(b_0)}$ modulo $p^d$ can
always be uniquely solved by choosing $y_1\equiv 0\pmod{p^d}$. Then,
for $b_0\in\{1,\cdots,p-1\}$, subtracting $a_0^{(b_0)}$ from
$f_{b_0}^*(y_1)$, one has
$f_{b_0}^*(y_1)-a_0^{(b_0)}=p\left(\sum_{k=d-1}^1p^{k-1}a_k^{(b_0)}y_1^k\right)=pf_{b_0}^{**}(y_1)$.
For $b_0=0$, one can make the similar operation to get
$f_{0}(x)-a_0^{(0)}=p\left(\sum_{k=d-1}^1p^{k-1}a_ky_1^k\right)=pf_0^{**}(y_1)$.
Apparently, $f_{b_0}(x)$ is uniquely determined by $a_0^{(b_0)}$
modulo $p^d$ and $f_{b_0}^{**}(y_1)$ modulo $p^{d-1}$. Applying the
hypothesis on $f_{b_0}^{**}(y_1)$, all equivalent polynomials of
$f_{b_0}^{**}(y_1)$ can be determined modulo $p^{d-1}$. Then, with
each valid\footnote{Note that not all equivalent polynomials of
$f_{b_0}^{**}(y_1)$ are valid, due to the existence of some power of
$p$ in each coefficient.} equivalent polynomial of
$f_{b_0}^{**}(y_1)$ and the value of $a_0^{(b_0)}$, one can further
uniquely determine all coefficients of $f(x)$ modulo $p^d$ in the
same way given in the proof of Lemma
\ref{lemma:determine-PP-dLessp}. Thus this theorem is proved.
\end{proof}

The above theorem tells us that all equivalent polynomials that
induce a given polynomial function modulo $p^d$ can be determined
via a recursive manner. Considering the complexity of solving Eq.
(\ref{equation:solve-PP-nLessp}) is $O(p^3)$ and the complexity of
solving Eq. (\ref{equation:solve-PP-pd}) is $O((pd)^3)$, the total
complexity of deriving one equivalent polynomial via the recursive
procedure is
\begin{equation}
O\left((pd)^3+(p(d-1))^3p+\cdots+\left(p\left(d-(d-p))\right)^3p^{d-p}\right)+p^3p^{d-p+1}\right)=O\left(p^{d-p+6}\right).
\label{equation:solve-PP-complexity1}
\end{equation}
In fact, this complexity can be further reduced, due to the
existence of a power of $p$ in each coefficient (except
$a_0^{(b_0)}$) of each sub-polynomial $f_{b_0}^*(y_1)$, which will
make more coefficients disappear as the value of $d$ decreases. In
the following, let us study what will happen when the sub-polynomial
and its derivatives are further decomposed.

At first, let us see the decomposition of the $p$ sub-polynomials
$\left\{f_{b_0}^*(y_1)=\sum_{k=d-1}^1p^{k-1}y_1^k\right\}_{b_0=0}^{p-1}$.
Similarly, choosing $y_1=py_2+b_1$, i.e., $y_2=\lfloor
y_1/p\rfloor=\lfloor x/p^2\rfloor$ and $b_1=(y_1\bmod
p)\in\{0,\cdots,p-1\}$, we have $p^2$ sub-polynomials as follows:
$\forall b_0,b_1\in\{0,\cdots,p-1\}$,
\[
f_{b_1,b_0}^*(y_2)=\begin{cases}%
\sum_{k=d-1}^1p^{2k-1}a_k^{(b_0)}y_2^k, & \mbox{when }b_1=0,\\
\sum_{k=d-1}^1\left(\sum_{l=d-1}^kp^{l-1}a_l^{(b_0)}\binom{l}{k}b_1^{l-k}p^ky_2^k\right)+\sum_{l=d-1}^1p^{l-1}a_l^{(b_0)},
& \mbox{when }b_1\in\{1,\cdots,p-1\}.
\end{cases}
\]
The above polynomial can be rewritten in the following form:
\[
f_{b_1,b_0}^*(y_2)=\begin{cases}%
\sum_{k=d-1}^1p^{2k-1}a_k^{(b_1,b_0)}y_2^k, & \mbox{when }b_1=0,\\
\sum_{k=d-1}^1p^{2k-1}a_k^{(b_1,b_0)}y_2^k+a_0^{(b_1,b_0)}, &
\mbox{when }b_1\in\{1,\cdots,p-1\},
\end{cases}
\]
where $a_k^{(0,b_0)}=a_k^{(b_0)}$ and
$a_k^{(b_1,b_0)}=\sum_{l=d-1}^kp^{l-k}a_l^{(b_0)}\binom{l}{k}b_1^{l-k}$
for $1\leq k\leq d-1$ and $1\leq b_1\leq p-1$. Then, solving
$a_0^{(b_1,b_0)}$ and subtracting it from the involved polynomial,
one can get $p^2$ polynomials modulo $p^{d-2}$ as follows:
\[
\left\{f_{b_1,b_0}^{**}(y_2)=\sum_{k=d-2}^1p^{2(k-1)}a_k^{(b_1,b_0)}y_2^k\right\}_{0\leq
b_0,b_1\leq p-1}.
\]

Repeat the above procedure for $i$ times, where $1\leq i\leq d-2$,
one can get $p^{i+1}$ polynomials modulo $p^{d-i}$: $\forall
b_0,\cdots,b_i\in\{0,\cdots,p-1\}$,
\[
f_{b_i,\cdots,b_0}^*(y_{i+1})=\begin{cases}%
\sum_{k=d-i}^1p^{i(k-1)+1}a_k^{(b_i,\cdots,b_0)}y_{i+1}^k, & \mbox{when }b_i=0,\\
\sum_{k=d-i}^1p^{i(k-1)+1}a_k^{(b_i,\cdots,b_0)}y_{i+1}^k+a_0^{(b_i,\cdots,b_0)},
& \mbox{when }b_i\in\{1,\cdots,p-1\},
\end{cases}
\]
and $p^{i+1}$ polynomials modulo $p^{d-i-1}$
\[
\left\{f_{b_i,\cdots,b_0}^{**}(y_{i+1})=\sum_{k=d-i-1}^1p^{i(k-1)}a_k^{(b_i,\cdots,b_0)}y_{i+1}^k\right\}_{0\leq
b_0,\cdots,b_i\leq p-1},
\]
where $y_{i+1}=\lfloor y_i/p\rfloor=\lfloor x/p^{i+1}\rfloor$ and
$b_i=(y_i\bmod p)\in\{0,\cdots,p-1\}$. Observing the above
polynomials, one can see that some higher coefficients disappear
modulo $p^{d-i}$ or modulo $p^{d-i-1}$ due to the existence of
powers of $p$. Assume $p^{i(k-1)+1}<p^{d-i}$ or
$p^{i(k-1)}<p^{d-i-1}$, one has $k<\frac{d-1}{i}$, so
$a_k^{(b_i,\cdots,b_0)}$ is valid only when
$k\leq\left\lceil\frac{d-1}{i}\right\rceil-1$. When
$i\geq\left\lceil\frac{d-1}{p}\right\rceil$, one has
$i\geq\frac{d-1}{p}\Rightarrow\frac{d-1}{i}\leq
p\Rightarrow\left\lceil\frac{d-1}{i}\right\rceil\leq p\Rightarrow
\left\lceil\frac{d-1}{i}\right\rceil-1\leq p-1$, so the coefficients
of $f_{b_i,\cdots,b_0}^*(y_{i+1})$ can be uniquely determined. This
means that the complexity of deriving one equivalent polynomial via
the procedure is reduced to be
\[
O\left((pd)^3+(p(d-1))^3p+\cdots+
\left(p\left(d-\left\lceil\frac{d-1}{p}\right\rceil\right)\right)^3p^{\left\lceil\frac{d-1}{p}\right\rceil}
+p^3p^{\left\lceil\frac{d-1}{p}\right\rceil+1}\right)\approx
O\left(d^3p^{\left\lceil\frac{d-1}{p}\right\rceil+3}\right),
\]
which is much smaller than Eq. (\ref{equation:solve-PP-complexity1})
when $d\gg p$.

\begin{remark}
From Theorem \ref{theorem:solving-poly-small-n}, solving the
coefficients of a polynomial of degree $\leq p$ needs only $p$
input-output values of the polynomial function. This means that it
is still possible to determine the polynomials if the induced
polynomial function is partially known. For the polynomials of
arbitrary degree modulo $p^d$, only
$p^{\left\lceil\frac{d-1}{p}\right\rceil+1}$ input-output values are
needed. Note that the needed input-output values should satisfy some
certain distribution modulo $p^d$, so the number of required
input-output values will be larger if the values are observed in a
random process\footnote{For example, if they are collected in a
known-plaintext attack to an encryption procedure based on a
permutation polynomials.}.
\end{remark}

\bibliographystyle{unsrt}
\bibliography{BPC}

\end{document}